\documentclass[11pt]{article}
\usepackage{amsmath,fullpage,latexsym,amssymb,version,graphicx,cite}
\usepackage{mathtools}
\usepackage{tikz,pgfplots}
\usetikzlibrary{shapes.multipart}
\usepackage{booktabs,multirow}
\usepackage[colorlinks=true]{hyperref}
\usepackage{newtxtext}
\usepackage[smallerops]{newtxmath}
\usepackage{enumitem}

\newcommand{\BEQ}{\begin{equation}}
\newcommand{\EEQ}{\end{equation}}
\newcommand{\BEAS}{\begin{eqnarray*}}
\newcommand{\EEAS}{\end{eqnarray*}}
\newcommand{\BEA}{\begin{eqnarray}}
\newcommand{\EEA}{\end{eqnarray}}
\newcommand{\BIT}{\begin{itemize}}
\newcommand{\EIT}{\end{itemize}}
\newcommand{\eg}{{\it e.g.}}
\newcommand{\ie}{{\it i.e.}}
\newtheorem{theorem}{Theorem}
\newtheorem{corollary}{Corollary}

\newcounter{algorithmctr}
\renewcommand{\thealgorithmctr}{\arabic{algorithmctr}}
\newenvironment{algdesc}[1]
  {\refstepcounter{algorithmctr}\begin{list}{}{%
      \setlength{\rightmargin}{.05\linewidth}%
      \setlength{\leftmargin}{.05\linewidth}}%
      \item[]{\setlength{\parskip}{0ex}\medskip\par
        \nopagebreak
        {\bf Algorithm \thealgorithmctr} \emph{(#1).}}}
  {{\setlength{\parskip}{-1ex}\nopagebreak\smallskip\par}\end{list}}

\newcommand{\reals}{\mathbb{R}}
\newcommand{\symm}{\mathbb{S}}
\newcommand{\SV}[2]{\symm^{#1}_{#2}}
\newcommand{\argmin}{\mathop{\rm arg min}}
\newcommand{\Rank}{\mathop{\rm rank}}
\newcommand{\dom}{\mathop{\rm dom}}
\newcommand{\intr}{\mathop{\rm int}}
\newcommand{\relint}{\mathop{\rm ri}}
\newcommand{\cl}{\mathop{\rm cl}}
\DeclareMathOperator{\lin}{\mathrm{lin}}
\DeclareMathOperator{\newspan}{\mathrm{span}}
\let\TeXspan\span
\makeatletter
\def\sp@n{\TeXspan\omit\advance\@multicnt\m@ne}
\makeatother
\renewcommand{\span}{\newspan}
\newcommand{\prox}{\mathrm{prox}}

\newcommand{\hp}{\phi_\mathrm{p}}
\newcommand{\dpr}{d_\mathrm{p}}
\newcommand{\hd}{\phi_\mathrm{d}}
\newcommand{\ddu}{d_\mathrm{d}}
\newcommand{\hpd}{\phi_\mathrm{pd}}
\newcommand{\dpd}{d_\mathrm{pd}}
\newcommand{\dd}{d}
\newcommand{\xavg}{x^\mathrm{avg}}
\newcommand{\zavg}{z^\mathrm{avg}}
\definecolor{mygreen}{RGB}{0,158,115}
\allowdisplaybreaks

\title{Bregman primal--dual methods for conic linear optimization via
exact-penalty and self-dual formulations}
\author{Xin Jiang%
\thanks{Department of Industrial and Systems Engineering,
University of Houston. Email: \texttt{xinjiang@uh.edu}}
\and Justin Shao%
\thanks{Ming Hsieh Department of Electrical and Computer Engineering,
University of Southern California. Email: \texttt{jmshao@usc.edu}}
\and Lieven Vandenberghe%
\thanks{Department of Electrical and Computer Engineering,
University of California, Los Angeles. Email: \texttt{vandenbe@ucla.edu}}
}
\date{November 15, 2026}

\begin{document}
\maketitle

\begin{abstract}
We study Bregman primal--dual proximal methods for conic linear 
programs (LPs).  For problems involving non-polyhedral cones, Bregman 
distances lead to generalized projections that are often simpler than 
Euclidean projections. 
In this paper, we reformulate conic LPs as equivalent exact-penalty problems
over truncated cones, and use self-dual embeddings to construct such 
formulations systematically.  
This approach yields Bregman proximal methods with simpler subproblems
and convergence guarantees over bounded domains.
We also study line searches for the Bregman primal--dual hybrid gradient
and extragradient methods applied to self-dual problems.
For logarithmically homogeneous kernels, the Bregman distances from
the starting points are unbounded on the truncated cone,
and the classical ergodic bounds do not apply.
For these kernels, we derive convergence rates for the
exact-penalty objective at the averaged iterates of both methods.
Numerical experiments with nonnegative trigonometric polynomial optimization
and sparse semidefinite optimization indicate the potential
of this framework as a general-purpose approach for incorporating
efficient Bregman projections in conic linear optimization.

\end{abstract}

\section{Introduction}
Primal--dual proximal methods are attractive for large-scale conic linear
optimization because they require only matrix--vector products
and projections on the primal or dual cones.
Bregman variants of these methods are especially appealing because they can
exploit problem geometry more effectively than Euclidean projections
and thereby reduce the complexity of the proximal updates.
However, many Bregman projections of interest are not naturally posed over
the entire cone; they are simpler when the cone is truncated, \ie,
intersected with a bounding halfspace.  
Motivated by this observation, we reformulate conic LPs as equivalent
exact-penalty problems over truncated cones, and use self-dual
embeddings to construct these bounded formulations systematically.
This yields Bregman primal--dual proximal methods with simpler subproblems
and formulations that better reflect the underlying conic structure.

We discuss primal and dual conic linear programs (LPs) in the form
\begin{subequations}
\label{e-conic-lp-pd}
\BEA 
\label{e-conic-lp-primal}
\mbox{Primal:} & & 
 \begin{array}[t]{ll}
 \mbox{minimize}   & c^T x \\
 \mbox{subject to} & Ax - b \in L \\
                   & x\in K,
 \end{array}  \\
\label{e-conic-lp-dual}
\mbox{Dual:} & & 
 \begin{array}[t]{ll}
 \mbox{maximize}   & b^T z \\ 
 \mbox{subject to} & -A^T z + c \in K^*  \\
                   & z \in L^*.
 \end{array}
\EEA
\end{subequations}
The variables are $x\in\reals^n$ and $z\in\reals^m$, 
and $K$ and $L$ are closed convex cones, not necessarily polyhedral 
or proper.  Two important special cases of the primal 
problem~(\ref{e-conic-lp-primal}) are LPs in standard form 
($K$ is a proper cone, $L =\{0\}$) 
and LPs in inequality form with free variables ($K = \reals^n$, 
$L$ is a proper cone).

\paragraph{Exact-penalty formulations}
Primal--dual proximal algorithms for the conic LPs require at each 
iteration matrix--vector products with $A$ and $A^T$, and projections on 
the cones $K$ and $L^*$.  They are well suited for large structured 
problems when these operations are 
inexpensive~\cite{ADHL+:26,AHLL:23,ADHL+:21,LXGY:25,LuY:25}.
To reduce the complexity per iteration further, 
variants that use non-Euclidean projections 
associated with generalized (Bregman) distances are of interest.  
Two fundamental types of primal--dual proximal algorithms are 
the \emph{primal--dual hybrid gradient method} (PDHG)
\cite{ChP:11,ChP:16} and the \emph{extragradient method} (EG)
\cite{Kor:76,Nem:04}, applied to the primal--dual optimality conditions or
the associated saddle-point problem.  For both classes of algorithms,
extensions that use Bregman distances have been developed 
\cite{JuN:12b,Tse:08,KoM:17,HaA:21,AuT:05}.

A convenient format for the study of proximal algorithms in
convex optimization is the primal--dual pair 
\begin{subequations}
\label{e-pd}
\BEA
\label{e-primal}
\mbox{Primal:} &  &
\begin{array}[t]{ll}
\mbox{minimize} & c^Tx + f(x) + g(b-Ax),
\end{array}  \\
\label{e-dual}
\mbox{Dual:} & &
\begin{array}[t]{ll}
\mbox{maximize} & b^Tz -g^*(z) - f^*(A^Tz-c),
\end{array}
\EEA
\end{subequations}
where $f:\reals^n \rightarrow \reals \cup \{+\infty\}$ and 
$g:\reals^m \rightarrow \reals \cup \{+\infty\}$
are closed convex functions, with conjugates $f^*$ and~$g^*$.
The primal and dual LPs~(\ref{e-conic-lp-pd}) are a special case 
of~(\ref{e-pd}) with
\BEQ \label{e-fg-lp}
f = \delta_K, \qquad g=\delta_{-L}, \qquad f^* = \delta_{-K^*}, 
\qquad g^* = \delta_{L^*},
\EEQ
where $\delta_C$ denotes the indicator function of the set $C$.
For this choice of $f$ and $g$, the primal--dual proximal algorithms 
require Euclidean or non-Euclidean projections on $K$ and $L^*$.
As an alternative to~(\ref{e-fg-lp}) one can apply the proximal algorithms
to an exact-penalty formulation of~(\ref{e-conic-lp-pd}).  Assume for
simplicity that $K$ and $L$ are proper cones.
Suppose strong duality holds for~(\ref{e-conic-lp-pd}), \ie,
the primal and dual optimal values are equal,
and that there exist primal and dual solutions $x^\star, z^\star$
that satisfy $\hat w^T x^\star < 1$ and $\hat s^T z^\star < 1$,
where $\hat w \in\intr(K^*)$ and $\hat s\in \intr(L)$. 
Then the primal and dual LPs can be shown to be equivalent to 
\begin{subequations}
\label{e-conic-lp-pd-ex-pen-2}
\BEA 
\label{e-conic-lp-pd-ex-pen-2-p}
\mbox{Primal:} &&
\begin{array}[t]{ll}
\mbox{minimize} & c^T x  + \delta_S(x) + \delta_T^*(b-Ax),
\end{array} \\
\label{e-conic-lp-pd-ex-pen-2-d}
\mbox{Dual:} &&
\begin{array}[t]{ll}
\mbox{maximize} & b^T z  - \delta_T(z) - \delta_S^*(A^Tz-c),
\end{array}
\EEA
\end{subequations}
where  $\delta_S$ and $\delta_T$ are the indicator functions of the truncated 
cones
\BEQ \label{e-ST-defs}
 S = \{x \in K \mid \hat w^Tx \leq 1\}, \qquad
 T = \{z \in L^* \mid \hat s^Tz \leq 1\}, 
\EEQ
and $\delta_S^*$, $\delta_T^*$ are their conjugates, \ie, the support 
functions of $S$ and $T$.  The last terms in the objective functions 
of~(\ref{e-conic-lp-pd-ex-pen-2-p}) and~(\ref{e-conic-lp-pd-ex-pen-2-d})
are exact penalties for the constraints $Ax-b \in L$ and
$-A^Tz + c\in K^*$ in the LPs~(\ref{e-conic-lp-primal}) 
and~(\ref{e-conic-lp-dual}).
The pair of problems~(\ref{e-conic-lp-pd-ex-pen-2}) is a special case 
of~(\ref{e-pd}) with 
\BEQ \label{e-fg-pen}
f= \delta_S, \qquad g = \delta_T^*, \qquad f^* = \delta_S^*, \qquad 
g^* = \delta_T.
\EEQ
The PDHG or EG algorithms applied to~(\ref{e-conic-lp-pd-ex-pen-2}) 
require at each iteration Euclidean or non-Euclidean projections on the 
sets $S$ and $T$.  
The choice~(\ref{e-fg-pen}) has advantages over~(\ref{e-fg-lp}) 
that are especially important for algorithms based on non-Euclidean 
distance functions.
\BIT
\item Generalized projections on the truncated cones $S$ and $T$ can be
easier to define and compute than projections on $K$ and~$L^*$.
When using generalized projections one must be able to guarantee 
that for every possible value of its arguments, 
the projection has a unique solution in the interior of the 
domain of the distance-generating function, as is the case for the 
squared Euclidean distance function and Euclidean projection on a closed 
convex set.  For a non-Euclidean projection this property depends
on the set and the generalized distance \cite{CeZ:97,BBC:03}, and is 
more difficult to satisfy when the set is unbounded.  
Specifically, for a conic LP with 
non-polyhedral cones $K$ and $L^*$, it may be difficult to find Bregman 
distances for which projections on the cones $K, L^*$ are uniquely
defined and easy to compute.  The task is simpler for the 
sets $S, T$ in~(\ref{e-ST-defs}).  The generalized projections discussed 
in \cite{BBT:17,Ced:25,ChV:18,JiV:22,LFN:18} illustrate this.  
Consider for example the distance generated by a logarithmically 
homogeneous barrier function for a proper cone~$K$.  
For this distance, projecting on $S$ requires solving an optimization 
problem
\[
 \begin{array}{ll}
 \mbox{minimize} & a^Tx + \phi(x) \\
 \mbox{subject to} & \hat w^T x \leq 1,
\end{array}
\]
which has a unique solution for every $a$. Without the inequality 
constraint, the projection only exists for $a\in\intr(K^*)$.

\item The domains $\dom f$ and $\dom g^*$ in~(\ref{e-fg-pen}) are bounded.
Explicit bounds on the feasible sets are an important assumption in the 
theory of first-order methods \cite[p.~125]{JuN:12a}.  For example, 
the convergence analysis of primal--dual first-order algorithms 
for~(\ref{e-pd}) typically concludes with the result that 
\BEQ \label{e-cvg-result-intro}
  \mathcal L(x_k,z) - \mathcal L(x,z_k)
  \leq \frac{\alpha}{k} (\|x-x_0\|_2^2 + \|z-z_0\|_2^2)
\quad \mbox{for all $x\in \dom f$, $z\in \dom g^*$},
\EEQ
where $\mathcal L = c^Tx + f(x)  + b^Tz - g^*(z) - z^TAx$ and $\alpha$ is
some constant (see, \eg, \cite{ChP:16} and Section~\ref{s-ls}).
If $\dom f$ and $\dom g^*$ are bounded, one can maximize the two sides of
the inequality over $x\in\dom f$ and $z\in\dom g^*$
to obtain an $O(1/k)$ convergence rate.
However, if the domains are not bounded, it is more difficult to
interpret~(\ref{e-cvg-result-intro}) as establishing convergence.

\item Most distance-generating functions for cones are strongly convex on
bounded sets of the form~(\ref{e-ST-defs}), but not on the entire cone.
This is important in convergence analyses of Bregman proximal methods
because strong convexity of the distance-generating function is a common 
assumption~\cite{ChP:16,Nem:04,AuT:05}.
For example, the negative entropy, which generates the relative entropy,
is strongly convex on the unit simplex but not on the nonnegative orthant.  

\item The PDHG and EG algorithms 
applied to~(\ref{e-conic-lp-pd-ex-pen-2}) 
generate primal and dual iterates $x_k, z_k$ that are feasible,
which simplifies stopping conditions.
On the other hand, when the algorithms are applied to~(\ref{e-pd}),
the primal iterates $x_k$ are in $K$ but may not satisfy the constraint
$Ax_k-b \in L$, and the dual iterates $z_k$ are in $L^*$
but may not satisfy $-A^Tz_k + c\in K^*$.
\EIT

\paragraph{Self-dual embedding}
The exact-penalty formulation~(\ref{e-conic-lp-pd-ex-pen-2}) requires 
vectors $\hat w \in \intr(K^*)$, $\hat s\in \intr(L)$ that satisfy
$\hat w^T x^\star < 1$ and $\hat s^T z^\star < 1$ for some 
optimal solutions $x^\star, z^\star$.
For some applications such bounds can be derived from problem structure.
In general, however, they are not available, and this limits the value of
the method as a general-purpose approach.  Moreover, overly loose bounds
(with small vectors $\hat w$, $\hat s$) may slow down the convergence or
cause numerical difficulties.  These difficulties can be overcome by
applying the exact-penalty formulation to a self-dual embedding of the
LPs~\eqref{e-conic-lp-pd}.  A self-dual conic LP has the form 
\BEQ \label{e-sd-lp}
\begin{array}{ll}
\mbox{minimize} & -b^T x \\
\mbox{subject to} & Ax - b \in K^* \\  & x \in K,
\end{array}
\EEQ
where $K$ is a closed convex cone and $A$ is a skew-symmetric matrix.
The self-dual property is evident if we note that~(\ref{e-sd-lp})
is a special case of~(\ref{e-conic-lp-pd}) with $b=-c$, $L=K^*$, 
and $A=-A^T$.  

Of the many interesting properties of self-dual LPs, 
the important one for our purpose is that the optimal set of 
a strictly consistent self-dual LP is nonempty and bounded,
and that bounds on the optimal set are readily available.  
Assume for simplicity that the cone
$K$ in~(\ref{e-sd-lp}) is proper (a self-dual LP over a cone that is not
proper can be written as an equivalent self-dual LP over a proper cone; see
Appendix~\ref{s-app-mixed-sd}), and that $\hat x$ is strictly
feasible, \ie, $\hat x\in \intr(K)$ and 
$\hat s = A\hat x-b \in\intr(K^*)$.
Then it can be shown that optimal solutions $x^\star$ satisfy 
\BEQ \label{e-sd-bounds}
\hat s^T x^\star  \leq \hat s^T\hat x.
\EEQ
As a consequence, an equivalent exact-penalty formulation of the self-dual 
LP is 
\BEQ \label{e-lp-sd-ex-penalty}
\begin{array}{ll}
\mbox{minimize} & -b^T x + \delta_C(x) + \delta_C^*(b-Ax) 
\end{array}
\EEQ
where $C = \{ x \in K \mid \hat s^T x \leq \mu\}$
and $\mu > \hat s^T \hat x$.
Problem~(\ref{e-lp-sd-ex-penalty}) is self-dual and an example of the
general self-dual convex problem 
\BEQ \label{e-prob-sd}
\begin{array}{ll}
\mbox{minimize} & -b^Tx + f(x) + f^*(b-Ax)
\end{array}
\EEQ
where $f$ is a closed convex function and $A$ is a skew-symmetric matrix.  

An important source of self-dual LPs is the \emph{self-dual embedding} of
the conic LPs~(\ref{e-conic-lp-pd}).  
In the Ye--Todd--Mizuno embedding \cite{YTM:94}
and its conic generalizations \cite{DRT:97,LSZ:96c},
a conic LP is transformed into an equivalent self-dual problem
with known strictly feasible points.  These points are constructed explicitly
from the data of the original problem and yield explicit bounds 
on the optimal set of the embedded problem.
Therefore, by applying the exact-penalty reformulation developed above 
to the embedded problem, one obtains an equivalent formulation 
over truncated cones.
This is the key property of the self-dual embedding used in this paper.

\paragraph{Contributions}
The main contribution of this paper is a systematic reformulation technique
for applying Bregman primal--dual proximal methods to conic linear programs.
We first show that, given suitable bounds, a conic LP can be replaced by
an equivalent exact-penalty formulation.
We then use self-dual embeddings to obtain the required bounds directly
from the problem data, avoiding problem-specific estimates of
primal and dual solution sets.
The resulting reformulations can be solved by Bregman proximal methods
and provide the bounded domains needed in the convergence analysis.
The second contribution is a convergence analysis of Bregman PDHG and
extragradient methods with line search, applied to self-dual problems.
For logarithmically homogeneous kernels, the Bregman distances from
the starting points are unbounded on the truncated cone, and the classical
ergodic bounds do not apply.
We show that the exact-penalty objective at the averaged iterates of
both methods converges at the rate $O((\log k)/k)$ for these kernels, and at
the rate $O(1/k)$ if the kernel is finite at a solution of the
exact-penalty problem.
Numerical results on nonnegative trigonometric polynomial optimization and
sparse semidefinite optimization suggest the potential of the self-dual
exact-penalty formulation as a general-purpose approach for cones with
efficient Bregman projections.

\paragraph{Related work}
Primal--dual proximal algorithms are widely used in convex optimization;
see, \eg, \cite{ChP:11,HeY:12,ChP:16,MaP:18}.  Their use in large-scale
linear programming has recently attracted renewed interest, driven in part
by the success of GPU implementations of PDHG-type methods 
\cite{ADHL+:21,AHLL:23,ADHL+:26,LuY:25}.
Recent work extends this direction to general conic linear optimization
\cite{LXGY:25,XF:26}.
These developments provide important motivation for this paper,
but our focus is different: we consider non-polyhedral conic linear programs
and use non-Euclidean proximal updates to exploit the structure of
the relevant conic constraints.

More specifically, this paper is related to Bregman primal--dual proximal
methods, which include the Bregman proximal point method, mirror-prox,
and Bregman PDHG
\cite{Eck:93,Gul:94,AuT:05,Nem:04,Tse:08,JuN:12a,JuN:12b,KoM:17,HaA:21}.
For conic linear optimization, generalized projections on truncated cones
can often be computed efficiently, even when the Euclidean projection on the
cone is expensive \cite{ChV:18,JiV:22,JiV:23,Ced:25}.
The truncation of cones requires bounds on the primal and dual optimal sets,
which in this paper are obtained from the problem data
by a self-dual embedding.

Self-dual optimization problems have been studied since the 1950s
\cite{GoT:56,Duf:56,Dor:61,Han:64}.  The paper by Duffin~\cite{Duf:56}
is an early reference on self-dual conic LPs with non-polyhedral cones.
Self-dual embeddings were introduced by Goldman and Tucker~\cite{GoT:56}
and later became a standard tool for initialization and
infeasibility detection in interior-point methods
\cite{YTM:94,JRT:94,XHY:96,DRT:97,LSZ:96c,Wri:97,Ye:97}.
They have also been used with operator splitting methods
and primal--dual proximal methods for conic optimization \cite{OCPB:16,Sha:24,LuT:26}.
The present paper uses self-dual embeddings for a different purpose.
Rather than using the embedding primarily to initialize an algorithm
or to detect infeasibility, we use it to produce the bounded formulations
needed by Bregman primal--dual proximal methods.

Adaptive stepsizes and line searches for primal--dual proximal methods
have been studied in both Euclidean and Bregman settings.
For Euclidean PDHG, such methods are proposed
in~\cite{GLYE+:15,MaP:18}, and PDLP uses a similar adaptive rule
heuristically~\cite{ADHL+:21,ADHL+:26}.
An extension of the Malitsky--Pock line search to a Bregman dual update
is developed in \cite{ChC:22}.
Yazdandoost Hamedani and Aybat \cite{HaA:21} analyze a line search for
Bregman PDHG with separate primal and dual stepsizes, while the methods
in \cite{JiV:22,JiV:23} combine a Bregman primal update with a Euclidean
dual update.  The line search for PDHG considered here uses a common kernel and stepsize
in the primal and dual updates of a self-dual problem.  For the
extragradient method we use the line search of
Tseng~\cite[section~5]{Tse:08}.  The relation between the line-search
techniques is discussed in Section~\ref{s-ls-alg}.

Barrier functions have also been used as distance-generating functions
or regularizers in the barrier subgradient method~\cite{Nes:11},
in online learning~\cite{AHR:12}, and in (proximal) mirror descent methods
\cite{DMT:26}.  The analysis in Section~\ref{s-ls-exact-pen} follows
\cite{Nes:11,AHR:12} in approximating boundary points by interior points.

\paragraph{Paper outline}

The rest of the paper is organized as follows.  Section~\ref{s-duality}
reviews duality and optimality conditions for the
problems~(\ref{e-conic-lp-pd}) and~(\ref{e-pd}), and their specialization
to self-dual problems.
Section~\ref{s-ex-pen} develops the bounded exact-penalty reformulation
for conic LPs and establishes its equivalence to the original primal--dual
pair.  Section~\ref{s-pd-prox} discusses Bregman primal--dual proximal
methods for the resulting formulations.
Section~\ref{s-ls} presents line searches for the Bregman PDHG and
extragradient methods on self-dual problems, and analyzes their convergence
for the exact-penalty formulation of a self-dual conic LP.
Section~\ref{s-numerical} reports numerical experiments,
and Section~\ref{s-conclusion} concludes the paper.
Appendix~\ref{s-app-mixed} extends the construction to mixed-form conic LPs,
Appendix~\ref{s-app-ytm} derives the corresponding formulation
from the Ye--Todd--Mizuno self-dual embedding, and Appendix~\ref{s-app-proofs}
contains the proofs of the convergence results in Section~\ref{s-ls}.

\section{Duality and optimality conditions} \label{s-duality}
In this section we first review the duality and optimality theory 
for the problems~(\ref{e-conic-lp-pd}) and~(\ref{e-pd}). 
These duality results are covered in \cite[section~31]{Roc:70} 
\cite[section~11.H]{RoW:98}.  We then apply the duality results to self-dual problems.

\subsection{Duality} \label{s-duality-rock}
The primal problem~(\ref{e-primal}) is called \emph{strongly consistent} if
\BEQ \label{e-strong-consistency}
 b \in \relint(A\dom f + \dom g) = A\relint(\dom f) + \relint(\dom g),
\EEQ
where $\relint$ stands for relative interior.  
The equality in~(\ref{e-strong-consistency}) follows from 
\cite[theorem 6.6]{Roc:70}.
The problem is \emph{strictly consistent} if it is strongly consistent
and the set $A\dom f + \dom g$ has nonempty interior, \ie, if
\BEQ \label{e-strict-consistency}
 b \in \intr(A\dom f + \dom g).
\EEQ
Note that $A\intr(\dom f) + \intr(\dom g) \subseteq 
\intr(A\dom f + \dom g)$ but equality does not hold in general.

\begin{theorem}%
[{\cite[corollary~31.2.1]{Roc:70}, \cite[theorem~11.39]{RoW:98}}] \hfill
\label{thm-dual}
\BIT
\item If problem~(\ref{e-primal}) is strongly consistent, then strong 
duality holds:
problems~(\ref{e-primal}) and~(\ref{e-dual}) have the same optimal value, 
which, for a strongly consistent problem, can be finite or $-\infty$.  
If problem~(\ref{e-primal}) is strongly consistent and the optimal value 
is finite, then the dual optimum is attained.

\item The optimal set of the dual problem~(\ref{e-dual}) is nonempty and 
bounded if and only if the optimal value is finite and the primal
problem~(\ref{e-primal}) is strictly consistent.
\EIT
\end{theorem}
Given the symmetry between~(\ref{e-primal}) and~(\ref{e-dual}), 
dual counterparts of these statements also hold.
For the dual problem, strong consistency means that
\[
c \in \relint(A^T \dom g^*  -\dom f^*) 
 = A^T \relint(\dom g^*) - \relint(\dom f^*)
\]
and strict consistency that $c \in \intr(A^T \dom g^*  -\dom f^*)$.

For the LP~(\ref{e-conic-lp-primal}), which is the special case with
$f = \delta_K$, $g= \delta_{-L}$, strict 
consistency~(\ref{e-strict-consistency}) requires
\BEQ \label{e-lp-strict-consistent}
  b \in \intr(AK - L).
\EEQ
This holds if and only if the following two conditions are satisfied.
\begin{enumerate}
\item The affine hull 
\[
\span(AK-L) = A\span(K) + \span(L)
\]
has dimension $m$ (the row dimension of $A$).
For a cone $C$, the affine hull is $C + (-C)$, the subspace spanned
by $C$, denoted here by $\span(C)$.  
The span of $C$ is also the smallest subspace containing $C$.
The subspace $A\span(K) + \span(L)$ has dimension $m$ if and only if
the orthogonal complement has dimension zero, \ie, if and only if
\BEQ \label{e-aff-hull-dim-m}
\left.\begin{array}{c}
y\in \lin(L^*)  \\ A^Ty \in \lin(K^*) 
\end{array}\right\}
\qquad \Longleftrightarrow \qquad 
y=0.
\EEQ
Here $\lin(C) = C\cap (-C)$ is the \emph{lineality space} of a 
convex cone $C$, \ie, the largest subspace contained in $C$.
We used the property that $\lin(C) = \span(C^*)^\perp$.

\item \emph{Strong consistency:}  
$b\in \relint(AK-L) = A\relint(K) - \relint(L)$.  This means there exists
an $\hat x$ that satisfies
\[
 \hat x \in \relint(K), \qquad A\hat x - b \in \relint(L).
\]
\end{enumerate}
For LPs in standard form or inequality form, strict consistency entails
strict feasibility and a rank condition.
Appendix~\ref{s-app-mixed} describes the details for a mixed-form LP,
with cones of the form $K = K_1 \times \reals^{n_2}$ and 
$L = L_1 \times \{0\}$ where $K_1, L_1$ are proper cones.

\subsection{Optimality conditions} \label{s-optimality-rock}
The algorithms discussed later in the paper are based on the following
characterization of primal and dual optimal solutions.
\begin{theorem}%
[{\cite[theorem 31.3]{Roc:70}, \cite[theorem 11.39]{RoW:98}}]
\label{thm-optimality}
Strong duality holds for the problems~(\ref{e-pd})
and $x, z$ are optimal for~(\ref{e-primal}) and~(\ref{e-dual}),
respectively, if and only if
\BEQ \label{e-pd-oc}
A^Tz - c \in \partial f(x), \qquad b-Ax \in \partial g^*(z).
\EEQ
\end{theorem}
Combined with Theorem~\ref{thm-dual} this theorem implies that if the 
primal problem~(\ref{e-primal}) is strongly consistent, then $x$ is an 
optimal solution if and only if there exists a $z$ that together with $x$ 
satisfies~(\ref{e-pd-oc}).

The optimality conditions~(\ref{e-pd-oc}) can be expressed in several 
equivalent ways, which form the basis of different types of algorithms.
\BIT
\item \emph{Monotone inclusion.}
The conditions~(\ref{e-pd-oc}) form a \emph{monotone inclusion}
\BEQ \label{e-mon-inc}
 0 \in F(y) + \partial h(y) 
\EEQ
where $y=(x,z)$,
\BEQ \label{e-prob-F-h-def}
 F(y) = 
\left[\begin{array}{cc} 0 & -A^T \\ A & 0 \end{array}\right]
 \left[\begin{array}{c} x \\ z\end{array}\right] + \left[\begin{array}{c}
   c \\ -b \end{array}\right], \qquad h(y) = f(x) + g^*(z),
\EEQ
and $\partial h(y) = \partial f(x) \times \partial g^*(z)$
is the subdifferential of the convex function $h$.
The mapping $F$ is monotone (affine with skew-symmetric linear term),
but not the gradient of a convex function (unless $A=0$).

\item \emph{Convex--concave saddle-point problem.}
The primal and dual problems in~(\ref{e-pd}) can be written as
\begin{subequations}
\label{e-prob-saddle}
\BEA 
\mbox{Primal:} & & 
\begin{array}[t]{ll}
 \mbox{minimize} & \sup\limits_{\tilde z}{\mathcal L(x,\tilde z)}
\end{array} \\
\mbox{Dual:} & &
\begin{array}[t]{ll}
 \mbox{maximize} & \inf\limits_{\tilde x}{\mathcal L(\tilde x,z)} 
\end{array}
\EEA
\end{subequations}
where $\mathcal L$ is the convex--concave function
\[
 \mathcal L(x,z) = c^T x + f(x) + b^Tz - g^*(z) - z^TAx.
\]
The inequality
\BEQ \label{e-saddle-weak}
 \sup_{\tilde z}{\mathcal L(x,\tilde z)}
 \geq \mathcal L(x,z) \geq
 \inf_{\tilde x}{\mathcal L(\tilde x,z)}
\EEQ
holds for all $x\in\dom f$ and $z\in \dom g^*$, without assumptions
on $f$ and $g$.  
Strong duality holds for the problems~(\ref{e-pd}) and $x, z$ are optimal
for~(\ref{e-primal}) and~(\ref{e-dual}), respectively, if and only if 
\BEQ \label{e-saddle}
  \sup_{\tilde z}{\mathcal L(x, \tilde z)} 
 = \mathcal L(x, z) = \inf_{\tilde x}{\mathcal L(\tilde x, z)}. 
\EEQ
This last condition is therefore equivalent to~(\ref{e-pd-oc}).

\item \emph{Generalized variational inequality.}
In a \emph{generalized variational inequality} problem one is given a 
closed convex function $h$ and a mapping $F$, and seeks 
a vector $y$ such that 
\BEQ \label{e-vi-gen}
h(\tilde y) - h(y) +  F(y)^T(\tilde y-y) \geq 0 
\quad \mbox{for all $\tilde y \in \dom h$}.
\EEQ
This generalizes a (standard) variational inequality problem,
in which $h$ is the indicator function of a convex set $C$,  
and~(\ref{e-vi-gen}) can be written as
\BEQ \label{e-vi}
F(y)^T(\tilde y-y) \geq 0 \quad \mbox{for all $\tilde y\in C$}.
\EEQ
The variational inequality problem is monotone if $F$ is a 
monotone mapping.

The condition~(\ref{e-vi-gen}) is the subgradient inequality 
that defines the subdifferential $\partial h(y)$ 
in~(\ref{e-mon-inc}).
This shows that~(\ref{e-pd-oc}) is equivalent to a generalized
variational inequality with $F$ and $h$ defined in~(\ref{e-prob-F-h-def}).
Using the skew-symmetric structure in $F$ the generalized
variational inequality problem can also be expressed as follows:
find $x,z$ such that
\[ 
\mathcal L(\tilde x,z)  - \mathcal L(x, \tilde z) 
\geq  0 \quad \mbox{for all $\tilde x\in \dom f$, $\tilde z\in \dom g^*$},
\]
\ie, $\inf_{\tilde x} {\mathcal L(\tilde x,z)} \geq
\sup_{\tilde z}{\mathcal L(x,\tilde z)}$. 
Since strict inequality would contradict~(\ref{e-saddle-weak}), this is
the same as~(\ref{e-saddle}).
\EIT

The optimality conditions for the LPs~(\ref{e-conic-lp-pd}) follow by
applying the general results to the cone indicator functions
$f=\delta_K$ and $g^*=\delta_{L^*}$.
The subdifferentials are $\partial f(x) = N_K(x)$ and 
$\partial g^*(z) = N_{L^*}(z)$, where $N_C$ denotes the 
\emph{normal cone} to a closed convex cone $C$, defined as
\[
 N_C(x) = \left\{ \begin{array}{ll}
 \emptyset  & x\not \in C \\
 \{y \in -C^* \mid x^Ty = 0 \} & x\in C.
 \end{array}\right.
\]
The optimality conditions~(\ref{e-pd-oc}) for the LPs~(\ref{e-conic-lp-pd}) 
can therefore be expressed as a conic \emph{linear complementarity problem}
\[
\left[\begin{array}{cc} 0 & -A^T \\ A & 0 \end{array}\right] 
\left[\begin{array}{c} x \\ z \end{array}\right]
+ \left[\begin{array}{c} c \\ -b \end{array}\right]
= \left[\begin{array}{c} w \\ s \end{array}\right],
\qquad
 x \in K, \quad z\in L^*, \quad s \in L, \quad w \in K^*, \qquad 
 w^Tx = s^Tz = 0.
\]
A convenient shorthand is
\BEQ \label{e-lcp}
K \ni x \perp (c-A^Tz) \in K^*, \qquad
L \ni (Ax-b) \perp z \in L^*.
\EEQ

\subsection{Self-dual problem} \label{s-sd-oc}
Problem~(\ref{e-primal}) is \emph{self-dual} if $A=-A^T$, $b=-c$, and
$g = f^*$. In that case the primal problem is formally the same as the 
dual~(\ref{e-dual}), except for the orientation of the objective.
The LP~(\ref{e-conic-lp-primal}) is self-dual 
if $A=-A^T$, $b=-c$, and $K,L$ are 
dual cones.

For a self-dual problem, the optimal value of the dual 
problem~(\ref{e-dual}) 
is the negative of the optimal value of the primal~(\ref{e-primal}).  
Therefore the optimal value is nonnegative (possibly $+\infty$ for
an infeasible problem) and strong duality holds if and only if 
the optimal value is zero.   We also note that if $x$ is primal optimal,
then one can take $z=x$ in the optimality conditions~(\ref{e-pd-oc}).
The following corollary summarizes the implications of 
Theorems~\ref{thm-dual} and~\ref{thm-optimality} for a self-dual problem.
\begin{corollary} Assume $A=-A^T$ in~(\ref{e-prob-sd}) and $f$ is
a closed convex function.
\BIT
\item
If problem~(\ref{e-prob-sd}) is strongly consistent, then the 
optimal value is zero and the optimum is attained.
\item The optimal set is nonempty and bounded if and 
only if the problem is strictly consistent.
\item The optimal value is zero and $x,\tilde x$ are optimal if and only
if
\[
 b-A x \in \partial f(\tilde x), \qquad b-A\tilde x \in \partial f(x). 
\]
\item The optimal value is zero and $x$ is optimal if and only if
\BEQ \label{e-sd-oc}
   b-A x \in \partial f(x).
\EEQ
\EIT
\end{corollary}
The optimality condition~(\ref{e-sd-oc}) is a monotone inclusion.
The other equivalent optimality statements listed in 
Section~\ref{s-optimality-rock} also simplify for a self-dual problem.
\BIT
\item \emph{Convex--concave saddle-point problem.}
For the self-dual problem, 
\BEQ \label{e-lag-sd}
 \mathcal L(x,\tilde x) = b^T(\tilde x -x) + f(x) - f(\tilde x) - 
 \tilde x^TAx.
\EEQ
Note that $\mathcal L(x,\tilde x) = -\mathcal L(\tilde x,x)$ and 
$\mathcal L(x,x) =0$.
The optimal value of~(\ref{e-prob-sd}) is zero and $x$ is optimal
if and only if 
\BEQ \label{e-sd-saddle}
  \sup_{\tilde x}{\mathcal L(x, \tilde x)} = 0
  = \inf_{\tilde x}{\mathcal L(\tilde x, x)}. 
\EEQ

\item \emph{Generalized variational inequality.}
For the self-dual problem, the generalized variational 
inequality~(\ref{e-vi-gen}) simplifies to the condition
\BEQ \label{e-sd-vi}
 f(\tilde x) - f(x) - (b-Ax)^T (\tilde x-x ) \geq 0 \quad
 \mbox{for all $\tilde x$}.
\EEQ
This is the subgradient inequality that defines the subdifferential
in~(\ref{e-sd-oc}).
The optimal value of~(\ref{e-prob-sd}) is zero and $x$ is optimal
if and only if~(\ref{e-sd-vi}) holds.
\EIT

An important example of a self-dual problem~(\ref{e-prob-sd}) is the 
LP~(\ref{e-sd-lp}) with $A=-A^T$.
This corresponds to $f = \delta_K$ and $f^* = \delta_{-K^*}$
in~(\ref{e-prob-sd}), and 
\BEQ \label{e-sd-lp-lag}
 \mathcal L(x,\tilde x) = b^T(\tilde x - x) - \tilde x^TAx + \delta_K(x) 
 - \delta_K(\tilde x).
\EEQ
The problem is strictly consistent if $b \in \intr{(AK - K^*)}$.
From the discussion in Section~\ref{s-duality-rock}, this holds if and only 
if the problem is strongly consistent, \ie, there exists an $\hat x$ that 
satisfies
\BEQ \label{e-sd-lp-strict-feas}
 \hat x \in \relint(K), \qquad  A \hat x - b\in \relint(K^*),
\EEQ
and the subspace 
\[
\span(AK-K^*) = A\span(K) + \span(K^*) 
= A \span(K) + \lin(K)^\perp
\]
has dimension $n$.  Equivalently, 
specializing~(\ref{e-aff-hull-dim-m}) to the
self-dual LP,
\BEQ \label{e-lin-K-K*}
 \left. \begin{array}{l}
 y\in \lin(K)  \\ Ay \in \lin(K^*) 
 \end{array}\right\} \qquad \Longleftrightarrow \qquad y =0.
\EEQ
A strictly consistent self-dual LP~(\ref{e-sd-lp}) satisfies the
following properties.
\BIT
\item In the optimality condition~(\ref{e-sd-oc}), $\partial f(x)$
is the normal cone $N_K(x) = \{ y \in -K^* \mid y^Tx = 0\}$ to $K$ 
at $x$. 
Hence $x$ is optimal if and only if $x$ and $s=Ax-b$ satisfy
the conic linear complementarity conditions
\[
 x \in K, \qquad s \in K^*, \qquad s^T x = 0.
\]

\item The optimal value is zero and the optimal set is nonempty and 
bounded.  Consider any $x$ and define $s = Ax-b$.
Let $\hat x$ be any point that satisfies~(\ref{e-sd-lp-strict-feas})
and define $\hat s = A\hat x-b$.
Taking the inner product with $\hat x -x$ on both sides of the 
equation $A (\hat x - x) = \hat s - s$
shows that $(\hat s - s)^T (\hat x- x) = 0$.  If $x$ is optimal, 
$s^T x= 0$ and 
\BEQ \label{e-normalizing}
 \hat s^T x + \hat x^T s = \hat s^T \hat x.
\EEQ
The right-hand side is positive because $\hat x\in\relint(K)$
and $\hat s \in\relint(K^*)$.

The fact that all optimal solutions satisfy~(\ref{e-normalizing})
gives a simple direct proof that the optimal set of~(\ref{e-sd-lp})
is bounded.  We note the following fact.  
If $C$ is a closed convex cone, then
\[
 \left.  \begin{array}{l}
 u\in C  \\ w \in \relint(C^*) \\ w^T u = 0  \end{array}
 \right\}
 \qquad \Longrightarrow \qquad u\in \lin(C).
\]
To show that the optimal set of~(\ref{e-sd-lp}) is bounded, suppose 
$x = x^\star + \alpha y$ is optimal for all $\alpha \geq 0$.
Then~(\ref{e-normalizing}) requires that 
\[
y\in K, \qquad Ay \in K^*, \qquad
\hat s^T y = 0, \qquad \hat x^T Ay = 0.
\]
The conditions $y\in K$, $\hat s^Ty =0$ imply $y\in\lin(K)$
because $\hat s \in\relint(K^*)$.
The conditions $Ay\in K^*$, $\hat x^TAy =0$ imply $Ay\in\lin(K^*)$
because $\hat x \in\relint(K)$.  By~(\ref{e-lin-K-K*}), $y$ is zero.
\EIT
Strict consistency and its implications are especially simple if 
$K$ is a proper cone.   In that case, the problem is strictly consistent if and 
only if there exists a strictly feasible $\hat x$, 
\ie, a point that satisfies 
$\hat x \in\intr(K)$ with $A\hat x- b\in\intr(K^*)$.
Note that properness of $K$ can be assumed without loss of generality.
When the cone is not proper, one may make its lineality space and span
explicit and reduce the problem to an equivalent self-dual LP 
over a proper cone; see Appendix~\ref{s-app-mixed}.

\section{Exact-penalty formulations} \label{s-ex-pen}
We now apply the results of the previous section to derive the exact-penalty
formulations~(\ref{e-conic-lp-pd-ex-pen-2}) and~(\ref{e-lp-sd-ex-penalty}).
Exact-penalty formulations have a long history, reviewed in~\cite{Bur:91}.
In the context of primal--dual proximal methods, they have appeared in 
convergence analyses and stopping criteria \cite{JiV:23,JiV:22,ShT:14}.
In this paper, they are used for equivalent problem formulations
to which the proximal algorithms can be directly applied.

\subsection{Minkowski gauges}
The functions $\delta_T^*$ and $\delta_S^*$ 
in~(\ref{e-conic-lp-pd-ex-pen-2}) are examples of \emph{Minkowski
gauges}.  The Minkowski gauge of a proper cone $C$ in $\reals^n$ is the 
function
\BEA
\sigma_C(y,a) 
& = & \inf{\{t\geq 0 \mid ta - y \in C\}} \label{e-minkowski-inf}\\
& = & \sup{\{y^Tu \mid a^Tu \leq 1, \, u\in C^*\}} \label{e-minkowski-sup}
\EEA
and is defined for all $a\in\intr(C)$ and $y\in\reals^n$.
The point $a$ is sometimes called the \emph{pole} of the Minkowski gauge.
The equivalence between the expressions~(\ref{e-minkowski-inf}) 
and~(\ref{e-minkowski-sup}) follows from conic duality.
The first expression~(\ref{e-minkowski-inf}) shows that 
$\sigma_C(y,a) = 0$ if $-y\in C$.
Otherwise, $1/\sigma_C(y,a)$ is the steplength from $a$ to the boundary
of the cone in the direction $-y$.
The second expression~(\ref{e-minkowski-sup}) shows that 
$\sigma_C(\cdot,a)$ is the support function of the set 
$\{u \in C^* \mid a^T u \leq 1\}$.
It is therefore also the conjugate of the indicator function of this set.
The Minkowski gauge $\sigma_C(\cdot, a)$ is a nonnegative, convex, and 
positively homogeneous function of its first argument, with domain 
$\reals^n$.

It is common in applications that $C$ is a direct product of proper cones.
Suppose $C = C_1 \times \cdots \times C_k$, where $C_1,\ldots, C_k$
are proper.  Then the Minkowski gauge for $C$ at
$y=(y_1, \ldots,y_k)$ and $a = (a_1, \ldots,a_k) \in \intr(C_1)
 \times \cdots \times \intr(C_k)$ is
\BEAS
\sigma_C(y,a)
 & = & \inf{\{ t \geq 0 \mid ta_i -y_i \in C_i, \; i=1,\ldots,k\}} \\
 & = & \sup{\{ \sum_{i=1}^k y_i^T u_i \mid \sum_{i=1}^k a_i^Tu_i \leq 1, 
 \; u_i \in C_i^*, \; i=1,\ldots,k\}}\\
 & = & \max{\{\sigma_{C_1}(y_1,a_1),  \ldots, \sigma_{C_k}(y_k,a_k)\}}.
\EEAS
We also note that if $C$ is a proper cone and $a\in \intr(C)$, then
\BEAS
\max{\{\sigma_C(y,a), \sigma_C(-y,a)\}}
 & = & \inf{\{ t \mid ta-y\in C,\, ta+y \in C\}} \\
 & = & \sup{\{ y^T (u-v) \mid a^T(u+v) = 1, \, u, v \in C^*\}}.
\EEAS
As a function of $y$, this is a norm.  If we denote the corresponding
unit norm ball by $B_a$, then $B_a+a$ is the largest convex set, symmetric
about $a$ and contained in $C$.  

As an example, for $C =\reals^m_+$ and a positive $m$-vector $a$, 
\[
 \sigma_C(y,a) = \max{\{0, \max_{i=1,\ldots,m} { y_i/a_i}\}},
\qquad
 \max{\{\sigma_C(y,a), \sigma_C(-y,a)\}}
   = \max_{i=1,\ldots,m}{ |y_i|/a_i}.
\]

\subsection{Exact penalties for conic linear optimization}
Consider the conic LPs~(\ref{e-conic-lp-pd}) with proper cones $K$ and $L$.
Assume that strong duality holds and that the primal and dual 
optima are attained.  Also suppose we are given points 
$\hat w\in\intr(K^*)$ and $\hat s \in \intr(L)$ that satisfy 
$\hat w^T x^\star < 1$ and $\hat s^T z^\star < 1$ for some primal 
optimal $x^\star$ and dual optimal $z^\star$.
Using the conjugate correspondence between Minkowski gauges and indicator
functions, we can write the problems~(\ref{e-conic-lp-pd-ex-pen-2})  as
\begin{subequations}
\label{e-lp-ex-pen}
\BEA
\label{e-lp-ex-pen-primal}
\mbox{Primal:} &&
\begin{array}[t]{ll}
\mbox{minimize} & 
c^T x + \sigma_{L}(b-Ax, \hat s) \\
\mbox{subject to} & x \in K,  \; \hat w^Tx \leq 1 
\end{array}  \\
\label{e-lp-ex-pen-dual}
\mbox{Dual:} &&
\begin{array}[t]{ll}
\mbox{maximize} & b^T z - \sigma_{K^*}(A^Tz -c, \hat w) \\
\mbox{subject to} & z\in L^*, \; \hat s^T z\leq 1.
\end{array}
\EEA
\end{subequations}
We will show that these problems have the same optimal value as
the primal and dual problems in~(\ref{e-conic-lp-pd}), respectively.
Moreover, every solution of~(\ref{e-lp-ex-pen-primal}) is a solution of
the primal LP in~(\ref{e-conic-lp-pd}), and every solution 
of~(\ref{e-lp-ex-pen-dual}) is a solution of the dual LP.  Conversely,
any optimal solution of the original primal or dual LP that satisfies the
corresponding imposed bound is optimal for the associated penalized problem.
The second term in the objective of~(\ref{e-lp-ex-pen-primal})
therefore acts as an exact penalty for the constraint $Ax-b\in L$, 
and the second term in the objective of~(\ref{e-lp-ex-pen-dual}) acts 
as an exact penalty for the constraint $-A^Tz+c \in K^*$ in the dual LP.

Using the definition~(\ref{e-minkowski-inf}), we can express the two 
problems~(\ref{e-lp-ex-pen}) as the pair of primal and dual conic 
LPs
\begin{subequations}
\label{e-conic-lp-pd-ex-pen}
\BEA 
\label{e-conic-lp-pd-ex-pen-primal}
\mbox{Primal:} &&
\begin{array}[t]{ll}
\mbox{minimize} & c^T x  + \alpha \\
\mbox{subject to} & Ax - b + \alpha \hat s \in L\\
 & x\in K, \; \hat w^T x \leq 1  \\
 & \alpha \geq 0
\end{array} \\
\label{e-conic-lp-pd-ex-pen-dual}
\mbox{Dual:} &&
\begin{array}[t]{ll}
\mbox{maximize} & b^T z  - \beta \\
\mbox{subject to} & -A^Tz + c + \beta \hat w \in K^* \\
     & z\in L^*, \; \hat s^T z \leq 1 \\
     & \beta \geq 0.
\end{array}
\EEA
\end{subequations}
In the primal problem we added a scalar variable $\alpha$ and in the dual
problem a scalar variable $\beta$. 
The equivalence between~(\ref{e-lp-ex-pen}) and~(\ref{e-conic-lp-pd-ex-pen}) 
follows directly from the definition of the Minkowski gauge.
The following observations establish the exactness of these reformulations,
(\ref{e-lp-ex-pen}) and~(\ref{e-conic-lp-pd-ex-pen}),
relative to the original conic LP pair~(\ref{e-conic-lp-pd}).
\BIT
\item The problems~(\ref{e-lp-ex-pen-primal}) 
and~(\ref{e-lp-ex-pen-dual}) are duals, and a special case 
of~(\ref{e-primal}) and~(\ref{e-dual}) with 
$f = \delta_S$, $g = \delta^*_T$, where $S$ and $T$
are defined in~(\ref{e-ST-defs}).
The support functions of these sets are
$f^* = \delta_S^* =\sigma_{K^*}(\cdot, \hat w)$
and $g = \delta_T^* = \sigma_L(\cdot, \hat s)$.

\item The problems~(\ref{e-lp-ex-pen-primal}) 
and~(\ref{e-lp-ex-pen-dual}) are strictly consistent.  
The primal problem~(\ref{e-lp-ex-pen-primal}) is strictly consistent
because $\dom\sigma_L(\cdot,\hat s) =\reals^m$ and there exist
points $x$ that satisfy $x\in\intr(K)$ and $\hat w^T x < 1$.
The dual problem~(\ref{e-lp-ex-pen-dual}) is strictly consistent
because $\dom\sigma_{K^*}(\cdot,\hat w) =\reals^n$ and there exist
points $z$ that satisfy $z\in\intr(L^*)$ and $\hat s^T z < 1$.
It follows from Theorem~\ref{thm-dual} that strong duality holds 
for~(\ref{e-lp-ex-pen-primal}) and~(\ref{e-lp-ex-pen-dual}), 
and that the optimum is attained in primal and dual.
Theorem~\ref{thm-optimality} provides necessary and sufficient
optimality conditions that are satisfied by these optimal solutions.

The same conclusions hold for the equivalent conic 
LPs~(\ref{e-conic-lp-pd-ex-pen}).
The optimality conditions for these problems are
\begin{subequations}
\label{e-lcp-penalty}
\BEQ
K \ni x \perp (\beta \hat w + c- A^Tz) \in K^*, \qquad
L \ni (\alpha \hat s + Ax-b) \perp z \in L^*
\EEQ
and
\BEQ
0 \leq \left[\begin{array}{c} \alpha  \\ 1- \hat w^Tx
\end{array}\right] \perp 
\left[\begin{array}{c} 1-\hat s^Tz \\ \beta \end{array}\right] \geq 0.
\EEQ
\end{subequations}

\item Suppose $x, z$ are optimal for~(\ref{e-conic-lp-pd}) and satisfy the
bounds $\hat w^Tx \leq 1$ and $\hat s^Tz \leq 1$.
Since $x$ is optimal for~(\ref{e-conic-lp-primal}),
it is feasible for~(\ref{e-lp-ex-pen-primal}) with
objective value $c^Tx + \sigma_L(b-Ax, \hat s) = c^T x$.
Since $z$ is optimal for~(\ref{e-conic-lp-dual}), it is feasible 
for~(\ref{e-lp-ex-pen-dual}) with objective value
$b^T z -\sigma_{K^*}(A^Tz -c,\hat w) = b^T z$.
By assumption, strong duality holds for~(\ref{e-conic-lp-pd}), so
$c^Tx= b^Tz$.   Therefore the two points $x,z$ are optimal 
for~(\ref{e-lp-ex-pen}), and the optimal value for~(\ref{e-lp-ex-pen})
is equal to the optimal value $c^Tx = b^Tz$ of~(\ref{e-conic-lp-pd}).

\item By assumption, there exists an optimal solution $x^\star$ 
of~(\ref{e-conic-lp-pd}) that satisfies $\hat w^T x^\star < 1$,  
and an optimal solution $z^\star$ of the dual problem
in~(\ref{e-conic-lp-pd}) that satisfies $\hat s^T z^\star < 1$.
Define $\alpha^\star = \beta^\star = 0$.
The optimality conditions~(\ref{e-lcp-penalty}) show that
$(x^\star,\alpha^\star)$ and $(z^\star, \beta^\star)$
are optimal solutions of~(\ref{e-conic-lp-pd-ex-pen}). 

Now suppose $x,z$ are any optimal solutions for~(\ref{e-lp-ex-pen}),
and define $\alpha = \sigma_L(b-Ax,\hat s)$ 
and $\beta = \sigma_{K^*} (A^Tz-c,\hat w)$. 
Since $(x,\alpha)$ is optimal for the primal 
problem~(\ref{e-conic-lp-pd-ex-pen-primal}), it must satisfy
the complementarity conditions~(\ref{e-lcp-penalty}) with 
$(z^\star,\beta^\star)$.  
Hence $\alpha=0$ and $x$ is feasible for~(\ref{e-conic-lp-pd}).
If $(z,\beta)$ is optimal for the dual 
problem~(\ref{e-conic-lp-pd-ex-pen-dual}), then it must
satisfy the complementarity conditions~(\ref{e-lcp-penalty}) with 
$(x^\star, \alpha^\star)$.  Hence $\beta=0$ and $z$ is feasible 
for~(\ref{e-conic-lp-pd}).  
Moreover $c^Tx = b^Tz$ because $x$ and $z$ are optimal 
for~(\ref{e-conic-lp-pd-ex-pen}) with $\alpha=\beta=0$ 
and strong duality holds for~(\ref{e-conic-lp-pd-ex-pen}).
Since $x$ and $z$ are feasible in~(\ref{e-conic-lp-pd}) and $c^Tx=b^Tz$
they are optimal for~(\ref{e-conic-lp-pd}).
\EIT

\subsection{Exact-penalty self-dual formulation}
We now turn to the self-dual LP~(\ref{e-sd-lp}).
Assume for simplicity that $K$ is a proper cone, 
and that $\hat x\in\intr(K)$ and $\hat s = A\hat x - b\in\intr(K^*)$.
Let $\mu > \hat s^T \hat x$.
Then~(\ref{e-normalizing}) shows that $\hat s^T x < \mu$ for all optimal
solutions $x$.  Therefore
\BEQ \label{e-sdlp-ex-pen}
 \begin{array}{ll}
 \mbox{minimize} & -b^T x + \mu \sigma_{K^*}(b-Ax, \hat s) \\
 \mbox{subject to} & x\in K, \; \hat s^T x \leq \mu,
 \end{array}
\EEQ
is an exact reformulation of the self-dual LP~(\ref{e-sd-lp}), \ie,
$x$ is an optimal solution of~(\ref{e-sdlp-ex-pen}) if and only if it is a 
solution of~(\ref{e-sd-lp}).
Note that~(\ref{e-sdlp-ex-pen}) is a self-dual problem of the
form~(\ref{e-prob-sd}) and can be written as
\[
\begin{array}{ll}
 \mbox{minimize}  & -b^T x + \delta_C(x) + \delta_C^*(b-Ax)
\end{array}
\]
where $C = \{x \in K \mid \hat s^T x \leq \mu\}$.
Problem~(\ref{e-sdlp-ex-pen}) is equivalent to the self-dual LP
\[
\begin{array}{ll}
\mbox{minimize} & -b^T x + \mu t \\ 
\mbox{subject to} &  \left[\begin{array}{cc}
      A & \hat s \\
    -\hat s^T & 0 \end{array}\right] \left[\begin{array}{c} x \\ t 
 \end{array}\right] - \left[\begin{array}{c} b \\ -\mu  \end{array} 
 \right] \in K^* \times \reals_+ \\
 & (x,t) \in K \times \reals_+.
\end{array}
\]
The objective function in~(\ref{e-sdlp-ex-pen}) can also be written as
\BEA
-b^T x + \mu \sigma_{K^*}(b-Ax, \hat s) 
+ \delta_C(x) 
& = & \sup_{\tilde x\in C}
 {(-b^T x + (b-Ax)^T \tilde x + \delta_C(x))} \nonumber \\ 
& = & 
\left\{ \begin{array}{ll}
\sup_{\hat s^T \tilde x\leq \mu} {\mathcal L(x,\tilde x)}
 & \hat s^T x \leq \mu \\
+\infty & \mbox{otherwise,} \end{array} \right.
\label{e-sdlp-gamma}
\EEA
where $\mathcal L$ is defined in~(\ref{e-sd-lp-lag}).
The identity~(\ref{e-sdlp-gamma}) will be used in
Section~\ref{s-ls-exact-pen} to show that the saddle-point gap
of~(\ref{e-sdlp-ex-pen}), restricted to $C$, is the sum of two
exact-penalty objective values.

\section{Primal--dual proximal methods} \label{s-pd-prox}
In this section we describe Bregman variants of the primal--dual
hybrid gradient (PDHG) method and the extragradient (EG) method,
and specialize them to the exact-penalty and self-dual formulations 
developed above.
The algorithms themselves are classical or closely related to existing
methods.

\subsection{Bregman PDHG and extragradient methods}
We discuss PDHG and EG applied to problems in the form~(\ref{e-pd}).
With fixed stepsizes, both algorithms can be interpreted as 
variations of the \emph{Bregman proximal point method} 
\cite{Gul:94,Eck:93,CeZ:97}, 
applied to the primal--dual optimality conditions~(\ref{e-mon-inc}).  Let 
\[
\dpd(y, \hat y) = \hpd(y) - \hpd(\hat y) - \nabla\hpd(\hat y)^T(y-\hat y)
\]
be a Bregman distance, generated by a differentiable convex function 
$\hpd: \reals^n\times \reals^m \rightarrow \reals$.  
In the Bregman proximal point method applied
to~(\ref{e-mon-inc}), the update $y_{k+1} = (x_{k+1}, z_{k+1})$ is 
defined as the solution of the inclusion problem
\BEA
\frac{1}{t} (\nabla \hpd(y_k) - \nabla \hpd(y_{k+1}))
 & \in & F(y_{k+1}) + \partial h(y_{k+1}) 
 \label{e-bregman-ppa} \\
 & = & \left[\begin{array}{c} c-A^Tz_{k+1} \\ Ax_{k+1}-b \end{array}\right]
    + \partial f(x_{k+1}) \times \partial g^*(z_{k+1}),
 \label{e-bregman-pd-prox}
\EEA
where $y_k = (x_k,z_k)$ and $t$ is a positive stepsize.
The position of $x_{k+1}$ and $z_{k+1}$ in the first term
on the right-hand side of~(\ref{e-bregman-pd-prox}) couples the two
blocks in the equation and complicates the computation of $y_{k+1}$.
In this section we discuss two classical approaches for simplifying the 
update.  The two methods will involve primal and dual Bregman distances 
\BEQ \label{e-d-pd}
 \dpr(x,x') = \hp(x) - \hp(x') - \nabla \hp(x')^T(x-x'),
\qquad
 \ddu(z,z') = \hd(z) - \hd(z') - \nabla \hd(z')^T(z-z'),
\EEQ
generated by convex functions $\hp :\reals^n\rightarrow \reals$ and 
$\hd:\reals^m\rightarrow \reals$, respectively.

\paragraph{Primal--dual hybrid gradient method}
The Bregman PDHG method can be interpreted as the Bregman proximal point
method with a special choice of Bregman distance 
\cite{HeY:12,JiV:22,JiV:23,ChP:16}. 
Consider the proximal point method with a distance-generating function
of the form
\BEQ \label{e-pdhg-kernel}
\hpd(x,z) = \hp(x) + \hd(z) \pm t\, z^TAx.
\EEQ
A condition for $\hpd$ to be convex is given in~(\ref{e-dpd-convex})
below.
The Bregman distance generated by~$\hpd$ is
\BEQ \label{e-pdhg-dist}
\dpd(y, y') = \dpr(x,x') + \ddu(z,z') \pm t\, (z-z')^T A(x-x').
\EEQ
If we choose the positive sign for $\pm$ in~(\ref{e-pdhg-dist}), 
substitution in~(\ref{e-bregman-pd-prox}) gives
\BEAS
 \tfrac{1}{t} (\nabla \hp(x_k) - \nabla\hp(x_{k+1}))
  + A^Tz_k-c & \in & \partial f(x_{k+1}) \\
 \tfrac{1}{t} (\nabla \hd(z_k) - \nabla\hd(z_{k+1}))
  + b - A(2x_{k+1}-x_k) & \in & \partial g^*(z_{k+1}).
\EEAS
This gives a first PDHG variant:
\begin{subequations}
\label{e-pdhg-ver1}
\BEA
x_{k+1} & = & \argmin_x 
  {\big((c-A^Tz_k)^Tx + f(x) + \tfrac{1}{t}\dpr(x,x_k) \big)} \\
z_{k+1} & = & \argmin_z {\big((A(2x_{k+1}-x_k) - b)^Tz + g^*(z) + 
 \tfrac{1}{t}\ddu(z,z_k) \big)}.
\EEA
\end{subequations}
The iteration~(\ref{e-pdhg-ver1}) consists of a primal update followed by 
a dual update,
and both updates require the evaluation of a generalized proximal operator.
With the negative sign in~(\ref{e-pdhg-dist}), we get
\BEAS
 \tfrac{1}{t} (\nabla \hp(x_k) - \nabla\hp(x_{k+1})) 
  + A^T(2z_{k+1}-z_k) -c  & \in & \partial f(x_{k+1}) \\
 \tfrac{1}{t} (\nabla \hd(z_k) - \nabla\hd(z_{k+1}))
  + b- Ax_k & \in & \partial g^*(z_{k+1}).
\EEAS
This is the second variant of PDHG: 
\begin{subequations}
\label{e-pdhg-ver2}
\BEA
z_{k+1} & = & 
  \argmin_z {\big((Ax_k-b)^Tz + g^*(z) + \tfrac{1}{t}\ddu(z,z_k) \big)} \\
x_{k+1} & = & \argmin_x {\big((c - A^T(2z_{k+1}-z_k))^Tx + f(x) + 
  \tfrac{1}{t}\dpr(x,x_k) \big)}.
\EEA
\end{subequations}
It is also the first variant applied to the dual problem~(\ref{e-dual}).
The algorithms~(\ref{e-pdhg-ver1}) and~(\ref{e-pdhg-ver2}) use the same
stepsize parameter $t$ in primal and dual updates.  Different 
primal and dual stepsizes can be accommodated via the definition of
$\dpr$ and $\ddu$.  For example, if quadratic kernels 
$\hp(x) = \tfrac{1}{2\tau}(x^Tx)$ and $\hd(z) = \tfrac{1}{2\sigma}(z^Tz)$
are used and $t=1$, the two versions of PDHG are
\BEAS
x_{k+1} & = & \prox_{\tau f} (x_k - \tau (c-A^Tz_k)) \\
z_{k+1} & = & \prox_{\sigma g^*} (z_k - \sigma (A(2x_{k+1}-x_k) - b))
\EEAS
and
\BEAS
z_{k+1} & = & \prox_{\sigma g^*} (z_k - \sigma (Ax_k-b))\\
x_{k+1} & = & \prox_{\tau f} (x_k - \tau (c - A^T(2z_{k+1}-z_k))).
\EEAS

If we apply the PDHG methods~(\ref{e-pdhg-ver1}) and~(\ref{e-pdhg-ver2})
to the exact-penalty conic LP formulation~(\ref{e-lp-ex-pen}), we obtain 
\BEAS
x_{k+1} & = & \argmin_{x \in K, \, \hat w^Tx\leq 1} 
 {\big((c-A^Tz_k)^Tx + \tfrac{1}{t}\dpr(x,x_k) \big)}\\
z_{k+1} & = & \argmin_{z\in L^*, \, \hat s^Tz \leq 1}
 {\big((A(2x_{k+1}-x_k) - b)^Tz + \tfrac{1}{t}\ddu(z,z_k) \big)}
\EEAS
and
\BEAS
z_{k+1} & = & \argmin_{z \in L^*, \, \hat s^Tz\leq 1}
 {\big((Ax_k-b)^Tz + \tfrac{1}{t}\ddu(z,z_k) \big)} \\
x_{k+1} & = & \argmin_{x\in K, \, \hat w^Tx \leq 1}
 {\big((c - A^T(2z_{k+1}-z_k))^Tx + \tfrac{1}{t}\dpr(x,x_k) \big)}.
\EEAS
These iterations require Bregman projections on truncations of
the cones $K$ and $L^*$. 

\paragraph{Extragradient method}
Korpelevich's extragradient method is a classical method for the monotone 
variational inequality problem~(\ref{e-vi}) \cite{Kor:76}, \cite[p.~285]{BeT:97a}.
Extensions of the extragradient method have been developed that
handle the generalized variational inequality~(\ref{e-vi-gen}) and use 
Bregman proximal operators.  The following algorithm was proposed and 
analyzed in \cite[section 3]{Nem:04}, \cite[section 4]{AuT:05},
\cite[section 5]{Tse:08}, \cite{JuN:12b,TN:26}.

The Bregman extragradient algorithm can be interpreted as two iterations 
of a fixed-point iteration to compute $y_{k+1}$ from~(\ref{e-bregman-ppa}).
We start by writing~\eqref{e-bregman-ppa} as a fixed-point condition 
for $y_{k+1}$:
\BEQ \label{e-eg-from-ppa}
 y_{k+1} = \argmin_y{ (h(y) + F(y_{k+1})^Ty + \tfrac{1}{t} \dpd(y,y_k))}.
\EEQ
Two steps in the fixed-point iteration to compute the fixed point $y_{k+1}$,
starting at $y_k$ as initial point, are
\begin{subequations}
\label{e-gen-eg}
\BEA
\bar y_k & = & \argmin_y{ (h(y) + F(y_k)^T y + \tfrac{1}{t} \dpd(y,y_k))}\\
 y_{k+1} & = & \argmin_y{ (h(y) + F(\bar y_k)^T y + \tfrac{1}{t} 
 \dpd(y,y_k))}.
\EEA
\end{subequations}
If we use a separable distance function $\dpd(y,y')=\dpr(x,x')+\ddu(z,z')$,
where $y=(x,z)$ and $y' = (x',z')$, and substitute the expressions for $F$
and $h$ in~(\ref{e-prob-F-h-def}), we obtain the algorithm%
\begin{subequations}
\label{e-extragradient}
\BEA
\bar x_k & = & \argmin_x {((c-A^Tz_k)^T x + f(x) 
  + \tfrac{1}{t} \dpr(x,x_k))} \\
\bar z_k & = & \argmin_z {((Ax_k-b)^T z + g^*(z)
  + \tfrac{1}{t} \ddu(z,z_k))} \\
x_{k+1} & = & \argmin_x {((c-A^T\bar z_k)^T x + f(x)
  + \tfrac{1}{t} \dpr(x,x_k))} \label{e-extragradient-c} \\
z_{k+1} & = & \argmin_z {((A\bar x_k - b)^T z + g^*(z)
  + \tfrac{1}{t} \ddu(z,z_k))}. \label{e-extragradient-d}
\EEA
\end{subequations}
With $t=1$ and quadratic kernels $\hp(x) = \frac{1}{2\tau} (x^Tx)$
and $\hd(z) = \frac{1}{2\sigma} (z^Tz)$, we get
\begin{subequations}
\label{e-extragradient-prox}
\BEA
\bar x_k & = & \prox_{\tau f} (x_k - \tau(c-A^Tz_k)) \\
\bar z_k & = & \prox_{\sigma g^*} (z_k - \sigma (Ax_k-b)) \\
x_{k+1} & = & \prox_{\tau f} (x_k - \tau(c-A^T \bar z_k)) \\
z_{k+1} & = & \prox_{\sigma g^*} (z_k - \sigma (A\bar x_k-b)).
\EEA
\end{subequations}
Applied to the exact-penalty conic LP formulation,
the extragradient method takes the following four steps:
\BEAS
\bar x_k & = & \argmin_{x\in K,\, \hat w^T x \leq 1}
     {((c-A^Tz_k)^T x + \tfrac{1}{t} \dpr(x,x_k))} \\
\bar z_k & = & \argmin_{z\in L^*,\, \hat s^T z \leq 1}
     {((Ax_k-b)^T z + \tfrac{1}{t} \ddu(z,z_k))} \\
x_{k+1} & = & \argmin_{x\in K, \, \hat w^T x \leq 1}
     {((c-A^T\bar z_k)^T x + \tfrac{1}{t} \dpr(x,x_k))} \\
z_{k+1} & = & \argmin_{z\in L^*, \, \hat s^T z \leq 1}
     {((A\bar x_k - b)^T z + \tfrac{1}{t} \ddu(z,z_k))}.
\EEAS
Each of these steps is a Bregman projection on a truncated cone.

\paragraph{Stepsize condition} \label{s-step}
In Section~\ref{s-ls} we will extend the two methods to incorporate adaptive
stepsizes and analyze their convergence on self-dual problems.
Here we mention standard results on allowable values of the fixed 
stepsize $t$ in~\eqref{e-pdhg-ver1},~\eqref{e-pdhg-ver2}, 
and~\eqref{e-extragradient}.

We assume that $\hp$ and $\hd$ are strongly convex with respect to norms
$\|\cdot\|_\mathrm p$ and $\|\cdot\|_\mathrm d$, normalized so that
\BEQ \label{e-distance-strong-cvx}
 \dpr(x,x') \geq \frac{1}{2}\|x-x'\|_\mathrm p^2, \qquad
 \ddu(z,z') \geq \frac{1}{2}\|z-z'\|_\mathrm d^2
\EEQ
for all $x,x'\in\dom f\cap\dom\hp$ with $x'\in\intr(\dom\hp)$, and
all $z,z'\in\dom g^*\cap\dom\hd$ with $z'\in\intr(\dom\hd)$.
The corresponding matrix norm is
\BEQ \label{e-A-norm}
 \|A\|
 = \sup_{x,y\neq 0}\frac{y^TAx}
 {\|x\|_\mathrm p\|y\|_\mathrm d}
 = \sup_{x\neq 0}\frac{\|Ax\|_{\mathrm d*}}{\|x\|_\mathrm p}
 = \sup_{y\neq 0}\frac{\|A^Ty\|_{\mathrm p*}}{\|y\|_\mathrm d},
\EEQ
where $\|\cdot\|_{\mathrm p*}$ and $\|\cdot \|_{\mathrm d*}$ are the
dual norms of $\|\cdot\|_\mathrm p$ and $\|\cdot\|_\mathrm d$.
In the two PDHG variants and the extragradient method, a fixed stepsize $t$
satisfying $t\|A\|\leq1$ can be used \cite{Nem:04,Tse:08,ChP:16,JiV:22}.
For PDHG, this condition guarantees convexity of the kernel~$\hpd$
in~(\ref{e-pdhg-kernel}), since~(\ref{e-distance-strong-cvx})
and~(\ref{e-A-norm}) imply
\BEQ \label{e-dpd-convex}
 \dpd(y,y') \geq
 \frac{1-t\|A\|}{2} (\|x-x'\|_\mathrm p^2+\|z-z'\|_\mathrm d^2)
 \geq 0.
\EEQ
For the extragradient method, $t\|A\|\leq 1$ guarantees the stepsize
condition used in the analysis of \cite[section~5]{Tse:08};
see \cite[lemma~1]{Tse:08}, which extends \cite[theorem~3.1]{Nem:04}.

\subsection{Specialization to self-dual problems} \label{s-pd-prox-sd}
For a self-dual problem ($b=-c$, $g=f^*$, and $A=-A^T$),
with equal primal and dual distances $\dpr = \ddu = \dd$,
the primal and dual variants of PDHG~(\ref{e-pdhg-ver1}) 
and~(\ref{e-pdhg-ver2}) are equivalent up to interchanging
the two sequences.  The first version is
\begin{subequations}
\label{e-pdhg-sd}
\BEA
x_{k+1} & = & 
\argmin_x{((Az_k-b)^Tx + f(x) + \tfrac{1}{t}\dd(x,x_k))}\\
z_{k+1} & = & \argmin_z{((A(2x_{k+1}-x_k) - b)^Tz + f(z) + 
   \tfrac{1}{t}\dd(z,z_k))}.
\EEA
\end{subequations}
This is identical to the second variant if we switch the notation
for the two sequences $x_k$ and $z_k$.
For the extragradient method~(\ref{e-extragradient}) applied to the
self-dual problem~(\ref{e-prob-sd}), if we choose $\dpr=\ddu = d$
and start the iteration~(\ref{e-extragradient}) at $x_0=z_0$, then 
$\bar x_k = \bar z_k$ and $x_k=z_k$ for all $k$, and the four steps
in the iteration collapse to two:\
\begin{subequations}
\label{e-extragradient-sd}
\BEA
\bar x_k & = & \argmin_x {((Ax_k-b)^T x + f(x) 
   + \tfrac{1}{t} d(x,x_k))} \label{e-extragradient-sd-1} \\
x_{k+1} & = & \argmin_x {((A\bar x_k - b)^T x + f(x)
   + \tfrac{1}{t} d(x,x_k))}.  \label{e-extragradient-sd-2} 
\EEA
\end{subequations}
This is also the extragradient method applied directly to the 
simplified optimality conditions~(\ref{e-sd-oc}).

For the exact-penalty formulation~(\ref{e-sdlp-ex-pen}) of the self-dual
conic LP~(\ref{e-sd-lp}), the PDHG and extragradient methods are
\begin{subequations}
\label{e-pdhg-sd-pen}
\BEA
x_{k+1} & = &
\argmin_{x\in K, \, \hat s^Tx \leq \mu}
 {((Az_k-b)^Tx + \tfrac{1}{t}\dd(x,x_k))}\\
z_{k+1} & = & \argmin_{z\in K, \hat s^Tz\leq \mu}
  {((A(2x_{k+1}-x_k) - b)^Tz + \tfrac{1}{t}\dd(z,z_k))},
\EEA
\end{subequations}
and
\begin{subequations}
\label{e-eg-sd-pen}
\BEA
\bar x_k & = & \argmin_{x\in K, \, \hat s^Tx\leq \mu}
  {((Ax_k-b)^T x + \tfrac{1}{t} d(x,x_k))} \\
x_{k+1} & = & \argmin_{x\in K, \, \hat s^Tx \leq \mu}
   {((A\bar x_k - b)^T x + \tfrac{1}{t} d(x,x_k))}.
\EEA
\end{subequations}

\section{Adaptive stepsizes for self-dual problems} \label{s-ls}
In this section we study adaptive stepsizes selected by a line search for PDHG
(\ref{e-pdhg-sd}) and the extragradient method~(\ref{e-extragradient-sd}),
applied to the self-dual problem~(\ref{e-prob-sd}).
The line search technique removes the need for an estimate of $\|A\|$
in the condition $t\|A\| \leq 1$ of Section~\ref{s-step}, and can potentially
allow stepsizes larger than $1/\|A\|$.
Throughout the section, the two methods are treated in parallel.
Section~\ref{s-ls-alg} presents a line search for each method.
The relation of the Bregman PDHG algorithm with line search
(Algorithm~\ref{alg-ls}) to existing line searches for Euclidean and 
Bregman PDHG is discussed after Algorithm~\ref{alg-ls}.
The line search for the Bregman extragradient method 
(Algorithm~\ref{alg-eg-ls}) was proposed in \cite[section~5]{Tse:08}, and is 
specialized here to the self-dual problem.

Section~\ref{s-ls-conv} presents ergodic convergence bounds for the two 
methods on bounded effective domains, following standard arguments for
Bregman proximal methods~\cite{Nem:04,Tse:08,ChP:16}.
The resulting bounds on the gap function require a finite Bregman radius,
\ie, the Bregman distances from the starting points must be bounded on
the effective domain.  
Section~\ref{s-ls-exact-pen} specializes the results to the exact-penalty
formulation of self-dual LPs.  
In this case the gap function reduces to the exact-penalty objective.
We remove the assumption of finite Bregman radius by allowing distances 
generated by logarithmically homogeneous functions, including
logarithmically homogeneous barriers.
For such distances, Theorem~\ref{thm-barrier} shows that the exact-penalty
objective at the averaged iterates converges to zero at the rate
$O((\log k)/k)$, even though the Bregman radius is infinite.
Theorem~\ref{thm-finite} shows that the rate improves to $O(1/k)$ if
the kernel is finite at a solution.  Both results hold for the extragradient 
method, and for PDHG if the stepsize increases are summable.

The analysis in this section is independent of the 
proximal point method interpretation in Section~\ref{s-pd-prox}.

\subsection{Bregman PDHG and EG with line search} \label{s-ls-alg}
We first consider the PDHG iteration~(\ref{e-pdhg-sd}) with
a variable stepsize $t_k$.

\begin{algdesc}{Bregman PDHG with line search} \label{alg-ls}
Choose starting points $x_0,z_0 \in\dom f\cap\intr(\dom\phi)$, 
parameters $\delta\in(0,1]$, $\beta\in(0,1)$, $\bar t>0$, 
an initial stepsize $t_{-1}\in(0,\bar t\,]$, 
and growth factors $\gamma_k \geq 1$, $k = 0,1,\ldots$.
At iteration $k$, set $t := \min{\{\bar t, \gamma_k t_{k-1}\}}$ and
repeat the following two steps.
\begin{enumerate}
\item Compute
\begin{subequations}
\label{e-ls-trial}
\BEA
x^+ & = & \argmin_x {((Az_k-b)^Tx + f(x) + \tfrac{1}{t}\dd(x,x_k))}
 \label{e-ls-trial-x} \\
z^+ & = & \argmin_z {((A(2x^+-x_k)-b)^Tz + f(z) + \tfrac{1}{t}\dd(z,z_k))}.
 \label{e-ls-trial-z}
\EEA
\end{subequations}
\item If
\BEQ \label{e-ls-test}
 -t\, (z^+-z_k)^T A (x^+-x_k) \leq
 \delta \left(\dd(x^+,x_k) + \dd(z^+,z_k)\right),
\EEQ
accept the stepsize and set $t_k := t$, $x_{k+1} := x^+$, and 
$z_{k+1} := z^+$.
Otherwise, set $t:=\beta t$ and return to step~1.
\end{enumerate}
\end{algdesc}
Each trial in Algorithm~\ref{alg-ls} requires the two generalized projections
in~(\ref{e-ls-trial}).  The acceptance test in step~2 requires one inner 
product and the evaluation of the Bregman distances 
$\dd(x^+,x_k)$ and $\dd(z^+,z_k)$.
Note that if $(z^+-z_k)^TA(x^+-x_k)$ is nonnegative, the trial stepsize
is accepted automatically.

For the squared Euclidean distance, the test coincides with
\cite[eq.~(29)]{GLYE+:15} for equal primal and dual stepsizes,
and is closely related to the adaptive rule used in PDLP \cite{ADHL+:21}.
Goldstein et al.\ \cite{GLYE+:15} establish an $O(1/k)$ ergodic rate
under a summability condition on the relative stepsize changes,
whereas the PDLP rule is used heuristically.
Malitsky and Pock \cite[eq.~(8)]{MaP:18} propose a different Euclidean
line search.  In the notation used here, their test compares
$\sqrt{\tau_k\sigma_k}\|A^T(z_{k+1}-z_k)\|$ with
$\|z_{k+1}-z_k\|$, and only the dual update is recomputed during backtracking.

For Bregman distances, the test is a pointwise version of the Chambolle--Pock 
condition \cite[eq.~(13)]{ChP:16}, specialized to equal stepsizes and 
problems with smooth terms; 
their condition is imposed for all pairs of points.
Yazdandoost Hamedani and Aybat \cite{HaA:21} consider general, possibly
nonlinear coupling functions and use Bregman distances in both updates.
In their merely convex setting, the accepted stepsizes are nonincreasing
and satisfy $\sigma_k=\gamma_0\tau_k$.
The methods in \cite[eq.~(36)]{JiV:22} and \cite[eq.~(49)]{JiV:23} instead use 
a Bregman distance in the primal update and a Euclidean distance in the dual
update.  The self-dual setting considered here uses the same distance
and stepsize in the two updates, so both updates can be evaluated 
by the same projection oracle.

An adaptive stepsize can be incorporated in the extragradient method 
\eqref{e-extragradient-sd} as follows. 
\begin{algdesc}{Bregman EG with line search} \label{alg-eg-ls}
Choose a starting point $x_0\in\dom f\cap\intr(\dom\phi)$, 
parameters $\beta\in(0,1)$, $\bar t>0$,
an initial stepsize $t_{-1}\in(0,\bar t\,]$, and 
growth factors $\gamma_k \geq 1$, $k = 0,1,\ldots$.
At iteration $k$, set $t := \min{\{\bar t, \gamma_k t_{k-1}\}}$ and
repeat the following two steps.
\begin{enumerate}
\item
Compute
\begin{subequations}
\label{e-extragradient-sd-alg}
\BEA
\bar x & = & \argmin_x {((Ax_k-b)^T x + f(x)
   + \tfrac{1}{t} d(x,x_k))} \label{e-extragradient-sd-alg-1} \\
x^+ & = & \argmin_x {((A\bar x - b)^T x + f(x)
   + \tfrac{1}{t} d(x,x_k))}.  \label{e-extragradient-sd-alg2}
\EEA
\end{subequations}
\item If
 \BEQ \label{e-tseng-sd}
 f(\bar x)-b^T\bar x \leq (A\bar x-b)^Tx^++f(x^+) +\frac{1}{t}\dd(x^+,x_k)
\EEQ
accept the stepsize and set $t_k:= t$, $\bar x_k := \bar x$, and
$x_{k+1} := x^+$.  Otherwise, set $t := \beta t$ and return to step 1.
\end{enumerate}
\end{algdesc}

\noindent
Each trial in Algorithm~\ref{alg-eg-ls} also requires two generalized
projections.  The acceptance test requires the function values $f(\bar x)$
and $f(x^+)$, one inner product, and the Bregman distance $\dd(x^+,x_k)$.
The test (\ref{e-tseng-sd}) is the condition \cite[eq.~(45)]{Tse:08}
specialized to~(\ref{e-extragradient-sd});
see also \cite[eq.~(3.4)]{Nem:04}.
Alternatively, Nemirovski~\cite[section~6]{Nem:04} adjusts the stepsize
adaptively in his numerical experiments.
Armijo-type line searches for Euclidean extragradient-type methods
are studied in \cite{Kho:87,IuS:97,Tse:00}.

\subsection{Convergence analysis} \label{s-ls-conv}
In this section we analyze the convergence of Algorithms~\ref{alg-ls} 
and~\ref{alg-eg-ls} for the self-dual problem~(\ref{e-prob-sd}). 
We make the following assumptions.
\begin{enumerate} [label=A\arabic*., ref = A\arabic*] 
\item \label{as-1}
$f$ is closed convex and $A = -A^T$.

\item \label{as-2}
 The set $D = \dom f \cap \dom \phi$ is bounded, and
 $D^\circ = \dom f\cap\intr(\dom\phi)$ is nonempty.

\item \label{as-3}
 The distance-generating kernel $\phi$ is strongly convex on $D$:
 there exists a norm $\|\cdot\|$ such that
 \BEQ \label{e-ls-str-cvx}
 \dd(u,v) \geq \frac{1}{2} \|u-v\|^2 \quad
 \mbox{for all $u \in D$, $v\in D^\circ$}.
 \EEQ

\item \label{as-4}
The minimization problems in~(\ref{e-extragradient-sd-alg})
and in~(\ref{e-ls-trial}) have unique solutions in $D^\circ$
for every stepsize $t>0$ and all $x_k, z_k \in D^\circ$.
\end{enumerate}
We define the constants
\BEQ \label{e-ls-rho}
 \rho_1 = \sup_{x\in D}{\dd(x,x_0)} + \sup_{z\in D}{\dd(z,z_0)}, \qquad
 \rho_2 = \sup_{u,v,p,q \in D}{|(u-v)^TA(p-q)|}.
\EEQ
We refer to the constant $\rho_1$ as the \emph{Bregman radius.}
The Bregman radius $\rho_1$ may be infinite.
By Assumption \ref{as-2}, the constant $\rho_2$ is finite: 
$\rho_2 \leq \|A\|_2 \,(\mathop{\rm diam} D)^2 < \infty$,
where $\|A\|_2$ is the spectral norm of $A$ and
$\mathop{\rm diam} D$ is the Euclidean diameter of $D$.

\paragraph{Termination of the line searches}
Let $\|A\|$ be the matrix norm induced by the norm in \ref{as-3}, 
\[
\|A\| 
 = \sup_{x,y\neq 0} \frac{y^TAx}{\|x\| \|y\|} 
 = \sup_{x\neq 0} \frac{\|Ax\|_*}{\|x\|}. 
\]
For the trial points $x^+$, $z^+$ in Algorithm~\ref{alg-ls},
\BEAS
 \left|(z^+-z_k)^TA(x^+-x_k)\right|
 & \leq & \|A\| \|z^+-z_k\| \|x^+-x_k\| \\
 & \leq & \frac{\|A\|}{2} \left(\|x^+-x_k\|^2 + \|z^+-z_k\|^2\right) \\
 & \leq & \|A\| \left(\dd(x^+,x_k) + \dd(z^+,z_k)\right)
\EEAS
by~(\ref{e-ls-str-cvx}).
Hence the test~(\ref{e-ls-test}) is satisfied whenever $t\leq\delta/\|A\|$.
Similarly, the test~(\ref{e-tseng-sd}) in Algorithm~\ref{alg-eg-ls}
is satisfied whenever $t\leq 1/\|A\|$ \cite[lemma~1]{Tse:08}. (The argument
in~\cite{Tse:08} only needs~(\ref{e-ls-str-cvx}) at the points $\bar x$,
$x^+$, and~$x_k$.)  
Since in each unsuccessful trial $t$ is multiplied by $\beta$,
only finitely many trial stepsizes are rejected.
The trial value at the start of iteration $k$ is at least
$t_{k-1}$, because $\gamma_k \geq 1$ and $t_{k-1}\leq\bar t$.
If the first trial is successful, $t_k\geq t_{k-1}$.
Otherwise, the last rejected value $t_k/\beta$ exceeds $\delta/\|A\|$,
so $t_k > \delta\beta/\|A\|$.
By induction, the accepted stepsizes in Algorithm~\ref{alg-ls} satisfy
\BEQ \label{e-ls-tmin}
 t_k \geq t_{\min} := \min{\{t_{-1}, \delta\beta/\|A\|\}},
 \quad k = 0,1,\ldots.
\EEQ
The same argument shows that the accepted stepsizes in
Algorithm~\ref{alg-eg-ls} satisfy~(\ref{e-ls-tmin}) with $\delta=1$.
In particular, if $\gamma_k=1$, $t_{-1}=\bar t=t\leq 1/\|A\|$ 
(and $\delta=1$ in Algorithm~\ref{alg-ls}), 
the first trial stepsize is accepted in every iteration, and the two 
algorithms reduce to the methods of Section~\ref{s-pd-prox-sd} with fixed 
stepsize~$t$.

\paragraph{Growth factors}
The factors $\gamma_k$ control the increase of the trial stepsize.  Define
\BEQ \label{e-ls-Gamma}
 \Gamma_k = \sum_{i=1}^k (t_i-t_{i-1})_+,
\EEQ
where $(a)_+ = \max\{a,0\}$.  Since $t_i\leq\min\{\bar t,\gamma_i t_{i-1}\}$,
we have $\Gamma_k\leq \bar t\sum_{i=1}^k(\gamma_i-1)$.
The following choices give useful bounds on $\Gamma_k$.
\begin{enumerate}[label=(\alph*), ref = \alph*] 
\item \label{g-a}
If $\gamma_k=1$, the stepsizes are nonincreasing and $\Gamma_k=0$.

\item \label{g-b}
If $\gamma_k=1+\alpha_k/t_{k-1}$, the first trial stepsize 
is $\min\{\bar t,t_{k-1}+\alpha_k\}$ and $\Gamma_k \le \sum_{i=1}^k\alpha_i$.
Hence $\Gamma_k$ is bounded if $\sum_k\alpha_k<\infty$.

\item \label{g-c}
If $\gamma_k=1+\alpha_k$, then 
$\Gamma_k \le \bar t \sum_{i=1}^k\alpha_i$.
Hence $\Gamma_k$ is bounded if $\sum_k\alpha_k<\infty$.

\item \label{g-d}
If $\gamma_k=1+\alpha(k+1)^{-p}$, where $\alpha>0$ and $p>0$, then 
\[
\Gamma_k= \left\{\begin{array}{ll}
  O(k^{1-p}) & 0<p<1 \\ 
  O(\log k)  & p=1 \\ 
  O(1) & p>1.
\end{array}\right.
\]
The choice $p=0.6$ is used in PDLP~\cite{ADHL+:21}.
\end{enumerate}

\paragraph{Convergence results}
For $u,v \in D$, define the gap function
\BEQ \label{e-ls-gap}
 \eta(u,v) = \sup_{z\in D}{\mathcal L(u,z)}
  - \inf_{x\in D}{\mathcal L(x,v)},
\EEQ
where $\mathcal L$ is the convex--concave function~\eqref{e-lag-sd}
associated with~(\ref{e-prob-sd}). Clearly, $\eta(u,v)\geq 0$
for $u,v\in D$, since $\mathcal L(u,u) = \mathcal L(v,v) = 0$.
Moreover, $\mathcal L(x,z) = -\mathcal L(z,x)$, so
\BEQ \label{e-ls-gap-diag}
 \eta(u,u) = 2 \sup_{z\in D}{\mathcal L(u,z)}.
\EEQ

The next theorem is the main convergence result for~Algorithm~\ref{alg-ls}.
For $k\geq 0$, we define the weighted averages of the iterates of
Algorithm~\ref{alg-ls} as
\BEQ \label{e-ls-avg}
 \xavg_k = \frac{1}{T_k}\sum_{i=0}^k t_i x_{i+1}, \qquad
 \zavg_k = \frac{1}{T_k}\sum_{i=0}^k t_i z_{i+1}, \qquad
 T_k = \sum_{i=0}^k t_i.
\EEQ
\begin{theorem} \label{thm-ls-ergodic}
Suppose Assumptions \ref{as-1}--\ref{as-4} hold.  
For all $x,z\in D$, the averaged 
iterates~(\ref{e-ls-avg}) of Algorithm~\ref{alg-ls} satisfy
\BEQ \label{e-ls-pointwise}
 \mathcal L(\xavg_k, z) - \mathcal L(x, \zavg_k) \leq
 \frac{\dd(x,x_0) + \dd(z,z_0) + 2\rho_2 (t_0 + \Gamma_k)}{T_k},
 \quad k = 0,1,\ldots.
\EEQ
If $\rho_1 < \infty$, then
\BEQ \label{e-ls-gap-bound}
 \eta(\xavg_k,\zavg_k) \leq \frac{\rho_1 + 2\rho_2(t_0+\Gamma_k)}{T_k},
 \quad k = 0,1,\ldots.
\EEQ
\end{theorem}
The proof is given in Appendix~\ref{s-thm-3}.
Since $T_k\geq(k+1)t_{\min}$, the right-hand side
of~(\ref{e-ls-gap-bound}) is $O(1/k)$ if $\Gamma_k = 0$
(case~(\ref{g-a})), and for the growth factors 
in~(\ref{g-b}) and~(\ref{g-c}) if $\sum_k\alpha_k < \infty$,
since then $\Gamma_k = O(1)$.
With the growth factors in~(\ref{g-d}), it is $O(k^{-p})$ for $0<p<1$,
$O((\log k)/k)$ for $p=1$, and $O(1/k)$ for $p>1$.
For the growth schedule of PDLP ($p=0.6$), it is $O(k^{-0.6})$.

For Algorithm~\ref{alg-eg-ls}, the averaged iterate is defined as
\BEQ \label{e-eg-avg}
 \xavg_k = \frac{1}{T_k}\sum_{i=0}^k t_i \bar x_i,
\EEQ
with $T_k$ defined in~(\ref{e-ls-avg}).
The next theorem is the ergodic bound of \cite[proposition~4(a)]{Tse:08}
(see also \cite[proposition~4.1]{AuT:05}),
specialized to the self-dual problem and written in a form that parallels
Theorem~\ref{thm-ls-ergodic}, using the identity
$\mathcal L(x,z)=-\mathcal L(z,x)$.
\begin{theorem}[{\cite[proposition~4(a)]{Tse:08}}] \label{thm-eg-ergodic}
Suppose Assumptions \ref{as-1}--\ref{as-4} hold.  
For $x,z\in D$, the averaged iterate~(\ref{e-eg-avg}) of
Algorithm~\ref{alg-eg-ls} satisfies
\BEQ \label{e-eg-sd-ergodic}
 \mathcal L(\xavg_k, z) - \mathcal L(x, \xavg_k) \leq
 \frac{\dd(x,x_0) + \dd(z,x_0)}{T_k},
 \quad k = 0,1,\ldots.
\EEQ
If $\rho_1 < \infty$, then
\BEQ \label{e-eg-sd-gap}
 \eta(\xavg_k, \xavg_k) \leq \frac{\rho_1}{T_k},
 \quad k = 0,1,\ldots.
\EEQ
In (\ref{e-eg-sd-gap}), $\rho_1$ is the constant in~(\ref{e-ls-rho})
with $z_0 = x_0$.
\end{theorem}
For completeness, a proof is given in Appendix~\ref{s-thm-4}. 
Since $T_k\geq(k+1)t_{\min}$, the right-hand side
of~(\ref{e-eg-sd-gap}) is $O(1/k)$.
The bound~(\ref{e-eg-sd-ergodic}) holds for any sequence of stepsizes
that satisfy~(\ref{e-tseng-sd}), and does not depend on how the stepsizes
change from one iteration to the next.
This is the main difference from Theorem~\ref{thm-ls-ergodic}.

The preceding results concern the weighted averages of the iterates and
do not establish convergence of the iterates.  Point convergence requires
additional assumptions on the kernel, such as boundedness of Bregman sublevel
sets and suitable continuity properties;
see, \eg, \cite[section~3.3.4]{JiV:22}.

\subsection{Exact-penalty formulation of a self-dual conic LP}
\label{s-ls-exact-pen}
We now specialize the results to the exact-penalty formulation
(\ref{e-sdlp-ex-pen}) of the self-dual LP~(\ref{e-sd-lp}),
\ie, to $f = \delta_C$ with
\BEQ \label{e-ls-C}
 C = \{x \in K \mid \hat s^Tx \leq \mu\}.
\EEQ
We make the following assumptions.
\begin{enumerate}[label=B\arabic*., ref = B\arabic*] 
\item \label{as-B1}
$A=-A^T$, $K$ is a proper cone, $\hat x\in\intr(K)$,
$\hat s = A\hat x - b \in\intr(K^*)$, and $\mu > \hat s^T\hat x$.

\item\label{as-B2}
The kernel $\phi$ is closed and convex, and $\intr(K)\subseteq\dom\phi$.

\item \label{as-B3}
The kernel $\phi$ is strongly convex on $C\cap\dom\phi$:
there exists a norm $\|\cdot\|$ such that
\[
 \dd(u,v) \geq \frac{1}{2}\|u-v\|^2 \quad
 \mbox{for all $u\in C\cap\dom\phi$, $v\in C\cap\intr(\dom\phi)$}.
\]

\item \label{as-B4}
For every $a\in\reals^n$, $t>0$, and
$y\in C\cap\intr(\dom\phi)$, the generalized projection
\[
 \argmin_{x\in C} {\big(a^Tx + \tfrac{1}{t}\dd(x,y)\big)}
\]
exists, is unique, and lies in $\intr(\dom\phi)$.
\end{enumerate}
Under \ref{as-B1}--\ref{as-B4}, 
Theorems~\ref{thm-ls-ergodic} and~\ref{thm-eg-ergodic}
apply to~(\ref{e-sdlp-ex-pen}) with $D = C\cap\dom\phi$.
The PDHG and extragradient methods for this problem are given
in~(\ref{e-pdhg-sd-pen}) and~(\ref{e-eg-sd-pen}).  Each update is a
generalized projection on the truncated cone~$C$.
Hence \ref{as-B3} is required only on the truncated cone.
This is important in applications that involve distances
that are strongly convex on the truncated cone but not on the entire cone,
as is the case for the examples in Section~\ref{s-numerical}.

For $f = \delta_C$, the function $\mathcal L$ is bilinear on
$C\times C$,
\BEQ \label{e-L-sdlp}
 \mathcal L(x,z) = b^T(z-x) - z^TAx
\EEQ
for $x, z\in C$, hence continuous.  By \ref{as-B2}, $\cl D=C$.  
The supremum and infimum in~(\ref{e-ls-gap}) are therefore
unchanged if $D$ is replaced by $C$.  We use the same notation for the
resulting extension to $C\times C$:
\[
 \eta(u,v) = \sup_{z\in C}{\mathcal L(u,z)}
  - \inf_{x\in C}{\mathcal L(x,v)}
\]
for $u,v\in C$.  By~(\ref{e-sdlp-gamma}) and the identity
$\mathcal L(x,z)=-\mathcal L(z,x)$,
\BEQ \label{e-ls-gap-exact-pen}
 \eta(u,v) = q(u) + q(v)
\EEQ
for $u,v\in C$, where $q(x) = -b^Tx + \mu \sigma_{K^*}(b-Ax,\hat s)$ is the 
objective function of the exact-penalty problem~(\ref{e-sdlp-ex-pen}).
Since $q(u)=\sup_{z\in C}\mathcal L(u,z)\geq\mathcal L(u,u)=0$,
both terms on the right-hand side of~\eqref{e-ls-gap-exact-pen}
are nonnegative.
The optimal value is zero, so $q(u)=0$ if and only if $u$ is optimal.

In the analyses of Euclidean PDHG \cite{GLYE+:15,ADHL+:21},
restricting the gap to a bounded set provides the finite diameter needed
in the convergence bound.  For the exact-penalty formulation,
the restricted gap also has a direct interpretation:
by~(\ref{e-ls-gap-exact-pen}),
$\eta(u,v)=0$ if and only if both $u$ and $v$ solve the self-dual LP
(\ref{e-sd-lp}).  This interpretation motivates analyzing the
exact-penalty formulation rather than~(\ref{e-sd-lp}).

\paragraph{Distances with finite Bregman radius}
Suppose \ref{as-B1}--\ref{as-B4} hold.
If $\rho_1<\infty$, Theorems~\ref{thm-ls-ergodic} and~\ref{thm-eg-ergodic},
combined with~(\ref{e-ls-gap-exact-pen}), apply directly.
The averaged iterates are in $D$, and $T_k\geq(k+1)t_{\min}$.  Therefore,
Algorithm~\ref{alg-ls} applied to~(\ref{e-sdlp-ex-pen}) satisfies
\BEQ \label{e-ls-exact-rate}
 q(\xavg_k) + q(\zavg_k) \leq
 \frac{\rho_1 + 2\rho_2(t_0 + \Gamma_k)}{(k+1)\,t_{\min}},
 \quad k = 0,1,\ldots.
\EEQ
Algorithm~\ref{alg-eg-ls} applied to~(\ref{e-sdlp-ex-pen}) satisfies
\BEQ \label{e-eg-exact-rate}
 2 q(\xavg_k) \leq \frac{\rho_1}{(k+1)\,t_{\min}}, \quad k = 0,1,\ldots,
\EEQ
with $z_0=x_0$ in the definition of $\rho_1$ and $\delta=1$ in $t_{\min}$.
In particular, the exact-penalty objective values at the averaged iterates
of Algorithm~\ref{alg-eg-ls} converge to zero at the rate $O(1/k)$.
For Algorithm~\ref{alg-ls}, they converge to zero if $\Gamma_k = o(k)$,
at the rates discussed after Theorem~\ref{thm-ls-ergodic}.

\paragraph{Distances generated by logarithmically homogeneous kernels}
The next theorem relaxes the condition $\rho_1<\infty$ under the
additional assumption that the distance
is generated by a \emph{logarithmically homogeneous} function.  
\begin{enumerate}[label=B\arabic*., ref = B\arabic*] 
\setcounter{enumi}{4}
\item \label{as-B5} $\phi$ is $\nu$-logarithmically homogeneous with respect
to $K$, where $\nu > 0$, \ie, 
\BEQ \label{e-log-hom}
\intr(K) \subseteq \dom \phi \subseteq K, \qquad
 \phi(\lambda x) = \phi(x) - \nu\log\lambda \quad
 \mbox{for all $x\in \intr(K)$ and $\lambda > 0$}.
\EEQ
\end{enumerate}
Examples include the logarithmically homogeneous barrier functions used in 
interior-point methods (in this case, $\dom \phi = \intr(K)$). 
The Itakura--Saito distance used in Section~\ref{s-num-ntp} is generated
by a logarithmically homogeneous function that is not a barrier function
(in this case $\dom \phi = K \setminus\{0\}$).

\begin{theorem} \label{thm-barrier}
Suppose Assumptions~\ref{as-B1}--\ref{as-B5} hold.
Define
\BEQ \label{e-rho3}
 \rho_3 = 2\phi(\hat x) - \phi(x_0) - \phi(z_0)
 + \sup_{u,v\in C}{\big(\nabla\phi(x_0)^T(x_0-u)
 + \nabla\phi(z_0)^T(z_0-v)\big)}.
\EEQ
Then Algorithm~\ref{alg-ls} applied to~(\ref{e-sdlp-ex-pen}) satisfies
\BEQ \label{e-barrier-rate}
 q(\xavg_k) + q(\zavg_k) \leq \frac{2\nu\log(k+2) + \rho_3 +
 2\rho_2(t_0+\Gamma_k)}{(k+1)\,t_{\min}}
 + \frac{2(\hat s^T\hat x + \rho_2)}{k+2}, \quad k = 0,1,\ldots,
\EEQ
and Algorithm~\ref{alg-eg-ls} applied to~(\ref{e-sdlp-ex-pen}) satisfies
\BEQ \label{e-barrier-rate-eg}
 2q(\xavg_k) \leq
 \frac{2\nu\log(k+2) + \rho_3}{(k+1)\,t_{\min}}
 + \frac{2(\hat s^T\hat x + \rho_2)}{k+2}, \quad k = 0,1,\ldots,
\EEQ
with $z_0 = x_0$ in the definition of $\rho_3$.
\end{theorem}
(Here the averaged iterates are defined in~(\ref{e-ls-avg})
and~(\ref{e-eg-avg}), and $t_{\min}$
in~(\ref{e-ls-tmin}), with $\delta=1$ for Algorithm~\ref{alg-eg-ls}.)
The constant $\rho_3$ is finite, because $C$ is bounded and
$\hat x, x_0, z_0\in\intr(\dom\phi)$.
It is nonnegative, because the choice $u=v=\hat x$ in the argument
of the supremum in~(\ref{e-rho3}) shows that 
$\rho_3 \geq \dd(\hat x,x_0) + \dd(\hat x,z_0)$.

For Algorithm~\ref{alg-ls}, the theorem shows that
$q(\xavg_k) + q(\zavg_k) = O((\nu\log k + \rho_2\Gamma_k)/k)$,
which is $O((\log k)/k)$ if $\Gamma_k = O(1)$.
For Algorithm~\ref{alg-eg-ls}, it shows that $q(\xavg_k) = O((\log k)/k)$.

The proof of Theorem~\ref{thm-barrier} is given in Appendix~\ref{s-thm-5}
and uses a technique from the analysis of a barrier subgradient method
\cite{Nes:11}, which was later applied to online learning with self-concordant
barrier regularizers~\cite{AHR:12}.
For mirror descent and proximal mirror descent with logarithmic barrier
kernels, convergence rates of order $O((\log k)/k)$ are derived
in a recent paper \cite{DMT:26}.

The next theorem establishes that the factor $\log(k+2)$ 
in~\eqref{e-barrier-rate} and~\eqref{e-barrier-rate-eg}
can be removed if $\phi$ is finite at a solution of~(\ref{e-sdlp-ex-pen}).
This is useful for applications involving logarithmically homogeneous
entropic kernels that are finite on most of the boundary of $K$. 
(See Section~\ref{s-num-ntp} for an example.)
\begin{theorem} \label{thm-finite}
Suppose Assumptions~\ref{as-B1}--\ref{as-B5} hold
and problem~(\ref{e-sdlp-ex-pen}) has a solution $x^\star$ with
$\phi(x^\star)<\infty$.  Define
\BEQ \label{e-rho4}
 \rho_4 = 2\phi(x^\star) - \phi(x_0) - \phi(z_0)
 + \sup_{u,v\in C}{\big(\nabla\phi(x_0)^T(x_0-u)
 + \nabla\phi(z_0)^T(z_0-v)\big)}.
\EEQ
Then Algorithm~\ref{alg-ls} applied to~(\ref{e-sdlp-ex-pen}) satisfies
\BEQ \label{e-finite-rate}
 q(\xavg_k) + q(\zavg_k) \leq
 \frac{2\left(2\nu\log 2 + \rho_4 + 2\rho_2(t_0+\Gamma_k)\right)}
 {(k+1)\,t_{\min}}, \quad k = 0,1,\ldots,
\EEQ
and Algorithm~\ref{alg-eg-ls} applied to~(\ref{e-sdlp-ex-pen}) satisfies
\BEQ \label{e-finite-rate-eg}
 2q(\xavg_k) \leq
 \frac{2\left(2\nu\log 2 + \rho_4\right)}{(k+1)\,t_{\min}},
 \quad k = 0,1,\ldots,
\EEQ
with $z_0 = x_0$ in the definition of $\rho_4$.
\end{theorem}
The constant $\rho_4$ is finite because $C$ is bounded and
$\phi(x^\star)<\infty$.  It is nonnegative, because the choice
$u=v=x^\star$ in the argument of the supremum in~(\ref{e-rho4}) shows that
$\rho_4\geq\dd(x^\star,x_0)+\dd(x^\star,z_0)$.
The right-hand side of~(\ref{e-finite-rate}) is $O(1/k)$ if $\Gamma_k=O(1)$.
The right-hand side of~(\ref{e-finite-rate-eg}) is $O(1/k)$. 

The proof of Theorem~\ref{thm-finite} is given in Appendix~\ref{s-thm-6}.

\section{Numerical experiments} \label{s-numerical}
We illustrate the reformulation and algorithmic framework on two classes of
conic linear programs.
Starting from a conic LP, we first form the self-dual embedding
described in Appendix~\ref{s-app-ytm}.
We then apply the exact-penalty reformulation to obtain a problem
over truncated cones.  Finally, we solve the resulting problem
with the Bregman primal--dual methods of
Sections~\ref{s-pd-prox} and~\ref{s-ls}.
The two examples differ only in the cone and in the Bregman projection
oracle: the first uses the cone of nonnegative trigonometric polynomials
with the Itakura--Saito distance, and the second uses
a sparse positive semidefinite cone with the distance generated
by a logarithmic barrier function.
The experiments are intended to demonstrate that the proposed
reformulation can be used as a practical route from a generic conic LP
to Bregman first-order methods when efficient generalized projections
are available.   

\subsection{Implementation details} \label{s-num-detail}
We describe the implementation for a conic LP in standard primal--dual form
\[
\begin{array}{llll}
  \mbox{minimize} & c^Tx & \mbox{maximize} & b^Tz \\
  \mbox{subject to} & Ax=b \qquad \qquad \qquad & 
    \mbox{subject to} & c-A^Tz \in K^* \\
  & x \in K,
\end{array}
\]
where $K$ is a proper cone.  The first step is to form the self-dual embedding
\begin{equation} \label{e-num-sd-embed}
\begin{array}{ll}
\mbox{minimize} & \theta \\[1mm]
\mbox{subject to} & \left[ \begin{array}{cccc}
  0 & -A^T & c & q_x \\
  A & 0 & -b & q_z \\
  -c^T & b^T & 0 & q_\tau \\
  -q_x^T & -q_z^T & -q_\tau & 0
\end{array} \right]
\left[ \begin{array}{c} x \\ z \\ \tau \\ \theta \end{array} \right]
+ \left[ \begin{array}{c} 0 \\ 0 \\ 0 \\ 1 \end{array} \right]
= \left[ \begin{array}{c} w \\ s \\ \kappa \\ \eta \end{array} \right] \\
& (x,z,\tau,\theta) \in K\times \reals^m \times \reals_+ \times \reals \\
& (w,s,\kappa,\eta) \in K^* \times \{0\}^m \times \reals_+ \times \{0\}.
\end{array}
\end{equation}
This is a self-dual conic LP in the form of \eqref{e-sd-lp}.
We then apply Bregman PDHG and extragradient methods to its exact-penalty
reformulation, as described in Section~\ref{s-pd-prox}.

The parameters in the embedding (\ref{e-num-sd-embed}) are defined as
\[
\left[ \begin{array}{c} q_x \\ q_z \\ q_\tau \end{array} \right]
= \frac{1}{\hat\theta} \left( \left[ \begin{array}{c}
  \hat w \\ 0 \\ 1
\end{array} \right] - \left[ \begin{array}{ccc}
  0 & -A^T & c \\ A & 0 & -b \\ -c^T & b^T & 0
\end{array} \right]
\left[ \begin{array}{c} \hat x \\ \hat z \\ 1 \end{array} \right] \right),
\]
where $\hat \theta = \hat w^T \hat x + 1$ and
$\hat x \in \intr(K)$, $\hat z \in \reals^m$, $\hat w \in \intr(K^*)$.
The points $\hat x, \hat z, \hat w$ are also used to initialize the Bregman 
proximal methods.
The choice of $\hat x, \hat z, \hat w$ in our experiments uses
the heuristic in CVXOPT~\cite{ADV:20}, which is adapted from 
the initialization strategy of~\cite{Meh:92},
and can be summarized as follows.
We assume that $A$ has full row rank and compute the least-norm primal estimate
$x_{\mathrm{ls}} = A^T (AA^T)^{-1} b$.
If this point is not in the interior of the cone,
it is shifted in a fixed interior direction:
\[
\hat x = \begin{cases}
  x_{\mathrm{ls}} & \alpha_{\mathrm p}<0 \\
  x_{\mathrm{ls}}+(1+\alpha_{\mathrm p})e_{\mathrm p} \quad 
	& \mbox{otherwise,}
\end{cases}
\]
where $e_\mathrm{p} \in \intr(K)$ and $\alpha_{\mathrm p} 
= \inf \{\alpha \mid x_{\mathrm{ls}}+\alpha e_{\mathrm p} \in K\}$.
Similarly, on the dual side, we compute
\[
z_{\mathrm{ls}} = \argmin_z \|c-A^Tz\|_2^2, \qquad
w_{\mathrm{ls}} = c-A^Tz_{\mathrm{ls}},
\]
and set $\hat z = z_\mathrm{ls}$ together with
\[
\hat w = \begin{cases}
  w_{\mathrm{ls}} & \alpha_{\mathrm d}<0 \\
  w_{\mathrm{ls}}+(1+\alpha_{\mathrm d})e_{\mathrm d} \quad &
	\mbox{otherwise,}
\end{cases}
\]
where $e_{\mathrm d} \in \intr(K^\ast)$ and $\alpha_{\mathrm d}
= \inf \{\alpha \mid w_{\mathrm{ls}} + \alpha e_{\mathrm d} \in K^*\}$.
The remaining components of the reference point are chosen as $\hat s=0$,
$\hat\tau = \hat\kappa = 1$, $\hat\theta = \hat w^T\hat x + 1$, 
and $\hat\eta = 0$.

The exact-penalty formulation of the embedded problem 
contains the bounded constraint
\begin{equation} \label{e-num-coupled-bnd}
  \hat w^T x+\tau \le \mu, \qquad x \in K, \qquad \tau\ge 0,
\end{equation}
where $\mu > \hat w^T\hat x + 1$.
Direct application of Bregman PDHG or extragradient methods couples 
the $x$- and $\tau$-updates through (\ref{e-num-coupled-bnd}).
In the implementation we replace (\ref{e-num-coupled-bnd})
by the weaker constraints
\[
  \hat w^Tx\le \mu, \qquad 0\le \tau\le \mu.
\]
Consequently, the $x$-update is separated from the $\tau$-update
and becomes the Bregman projection
\begin{equation} \label{e-num-proj}
\Pi(a,\hat w,\mu) = \argmin_{\hat w^Tx \leq \mu, \, x\in K}
{(-a^Tx+\phi(x))},
\end{equation}
where $\phi$ is the distance-generating function.

PDHG uses the line search of Algorithm~\ref{alg-ls} 
and the extragradient method uses
Algorithm~\ref{alg-eg-ls}.  In both methods the trial value at the start of
iteration $k$ is $t = \min\{\bar t, \gamma_k t_{k-1}\}$,
with $\bar t = t_{-1} = 1$ and $\gamma_k = 1 + (k+1)^{-1.1}$.
The backtracking factor is $\beta = 0.8$, and the parameter
in the line-search test~(\ref{e-ls-test}) is $\delta = 1$.

At each iteration with $\tau>0$, we form the recovered primal--dual
variables as described at the end of Appendix \ref{s-app-ytm-embed}:
\[
\bar x = x/\tau, \qquad \bar z = z/\tau, \qquad \bar w = w/\tau.
\]
The residuals are measured for the original conic LP:
\[
r_{\mathrm p} = \frac{\|A\bar x-b\|_2}{1+\|b\|_2}, \qquad
r_{\mathrm d} = \frac{\|A^T\bar z+\bar w-c\|_2}{1+\|c\|_2}, \qquad
r_{\mathrm{gap}} = \frac{|c^T\bar x-b^T\bar z|}{1+|c^T\bar x|+|b^T\bar z|}.
\]
The algorithm is terminated when $r_\mathrm{max} =
\max\{r_{\mathrm p}, r_{\mathrm d}, r_{\mathrm{gap}}\} \leq 10^{-3}$
or when the number of iterations reaches $10^6$.  

All experiments are carried out 
in Python 3.14 on a Mac with an Apple M2 chip and 16\,GB of RAM.

\subsection{Nonnegative trigonometric polynomials} \label{s-num-ntp}
We first consider conic LPs over the cone of nonnegative trigonometric
polynomials of degree at most~$p$:
\[
K = \{x \in \reals^{p+1} \mid F_x(\omega) \ge 0 \mbox{ for all } \omega \},
\]
where $F_x(\omega) = x_1 + 2\sum_{k=2}^{p+1} x_k \cos ((k-1) \omega)$.
Optimization problems over the cone of nonnegative trigonometric polynomials
or its dual cone, the cone of positive semidefinite Toeplitz matrices,
are common in signal processing and system identification.
As discussed in \cite{ChV:18}, a useful Bregman distance on $K$ is
the Itakura--Saito distance generated by the negative Burg entropy
\[
\phi(x) = -\frac{1}{2\pi} \int_0^{2\pi} \log F_x(\omega) \, d\omega.
\]
In \cite{ChV:18}, an efficient method is described to compute a generalized
projection \eqref{e-num-proj}.
The key observation is that the optimality conditions for this projection
reduce to a one-dimensional nonlinear equation \cite[section~5]{ChV:18}.
Each Newton step for this scalar nonlinear equation requires
the solution of a positive definite Toeplitz linear system, which can be
carried out with complexity $O(p^2)$ by the Levinson--Durbin algorithm,
or with complexity $O(p(\log p)^2)$ using superfast Toeplitz solvers.
Hence the cost of evaluating the Bregman projection is roughly quadratic
in the degree $p$ of the trigonometric polynomial,
with a very small number of Newton steps in practice.
This should be contrasted with the Euclidean projection on the same
truncated cone, which is a significantly more difficult convex optimization
problem and, through its semidefinite representation, is tied to eigenvalue
decomposition or SDP subproblems of order $p+1$.

The Itakura--Saito kernel is $1$-logarithmically homogeneous
with $\dom\phi = K\setminus\{0\}$, \ie, $\phi$ is finite
on the boundary of $K$, except at the origin.
Theorem~\ref{thm-finite} therefore applies if the exact-penalty problem
has a nonzero solution.  With the growth factors in
Section~\ref{s-num-detail}, $\Gamma_k = O(1)$, so the exact-penalty
objective converges at the rate $O(1/k)$ for both methods.

\paragraph{Numerical results}
Table~\ref{t-ntp}%
\begin{table}
\centering
\begin{tabular}{ccrrrrr} \toprule
  & & & &
  &
  \multicolumn{2}{c}{Projection statistics} \\
  \cmidrule(lr){6-7}
  $p$ & $m$ & \multicolumn{1}{c}{Method} & \# iter.\ &
  Time/iter.\ (s) &
  active & Newton \\
	\midrule
  \multirow{2}{*}{$1000$} & \multirow{2}{*}{$200$}
	& EG   &  7553 & 0.01325 & 7 & 5.71 \\
	& & PDHG & 5728 & 0.01931 & 8 & 7.89 \\
	\midrule
  \multirow{2}{*}{$2000$} & \multirow{2}{*}{$500$}
	& EG   &  4710 & 0.03573 & 12 & 5.25 \\
	& & PDHG & 3957 & 0.03817 & 21 & 4.91 \\
	\midrule
  \multirow{2}{*}{$5000$} & \multirow{2}{*}{$1000$}
	& EG   & 10918 & 0.2164  & 6 & 5.83 \\
	& & PDHG & 11948 & 0.2917 & 11 & 5.76 \\
	\midrule
  \multirow{2}{*}{$6000$} & \multirow{2}{*}{$1200$}
	& EG   &  9759 & 0.4967 & 5 & 6.20 \\
	& & PDHG & 7562 & 0.5194 & 6 & 5.98 \\
	\bottomrule
\end{tabular}
\caption{Results for randomly generated LPs over the cone of nonnegative
trigonometric polynomials.  $p$ is the degree and $m$ is the number of
equality constraints.  The column \emph{active} gives the number of
projections for which $\hat w^Tx=\mu$ at the solution (the total number
of projections is twice \# iter.).  \emph{Newton} gives the average number
of Newton steps (averaged over the active projections).
Inactive projections are available in closed form.}
\label{t-ntp}
\end{table}
reports results for randomly generated linear programs over the cone of
nonnegative trigonometric polynomials.  The table shows that the per-iteration 
complexity remains well controlled as the problem size increases.
The bound $\hat w^Tx\leq\mu$ is active at only a small fraction of the
projections, and at those projections five to eight Newton steps suffice,
with no apparent dependence on $p$.  The truncation introduced by the
exact-penalty reformulation is therefore essentially free: it bounds the
domain, as the convergence theory requires, without materially changing the
cost of the proximal step.  The time per iteration grows roughly
quadratically with $p$, consistent with the cost of the projection routine
for the cone of nonnegative trigonometric polynomials \cite{ChV:18}.

\begin{figure}
\centering
\includegraphics{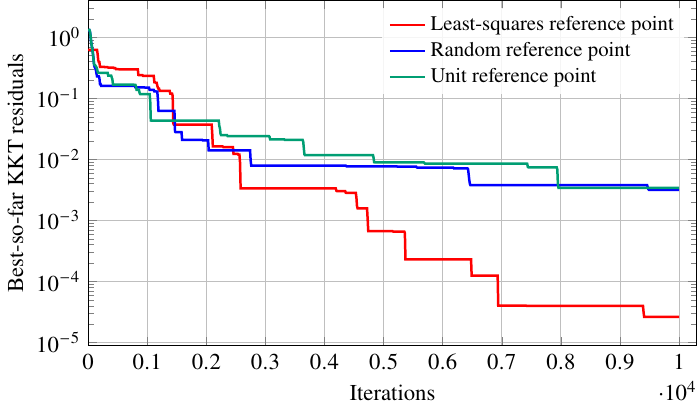}
\caption{Best-so-far KKT residuals for the Bregman extragradient method,
using three different reference points $(\hat x,\hat z,\hat w)$ in the
self-dual embedding~(\ref{e-num-sd-embed}).}
\label{f-ntp-inits}
\end{figure}

Figure~\ref{f-ntp-inits}%
compares three choices of the reference point $(\hat x, \hat z,
\hat w)$ on a representative instance, using the Bregman extragradient
method with line search.  The reference point enters the
embedding~(\ref{e-num-sd-embed}) via the definition of
$q_x, q_z, q_\tau$, and is also used to initialize the algorithm.
The first choice is the \emph{least-squares reference point} described in
Section~\ref{s-num-detail}, obtained from the least-norm primal and
least-squares dual estimates and shifted into the cone.  
The second choice is a \emph{random reference point}.
The third is the \emph{unit reference point}
\[
\hat x = (1,0,\ldots,0) \in \intr(K), \qquad
\hat w = (1,0,\ldots,0) \in \intr(K^*), \qquad
\hat s = 0,
\]
with a randomly chosen $\hat z$.  The vector $(1,0,\ldots,0)$ 
corresponds to the constant polynomial $F_x(\omega)=1$.
The figure shows the best-so-far KKT residuals versus iteration number.
Although the least-squares reference point is only a heuristic, it gives
substantially better performance in this example and is also more reliable
in most of our experiments.  We therefore use it in all the tests.

\subsection{Sparse semidefinite programming} \label{s-num-sdp}
We next consider linear programs over sparse semidefinite cones
\cite{ADV:13,VdA:15,JiV:22}. 
Let $E$ be an $n\times n$ chordal sparsity pattern, and let $\SV{n}{E}$
denote the subspace of symmetric matrices with zeros outside $E$.
We consider the primal--dual pair of cones
\[
K = \Pi_E(\mathbb S^n_+), \qquad K^* = \SV{n}{+} \cap \SV{n}{E},
\]
where $\Pi_E$ denotes projection onto the entries in $E$.
The primal cone $K$ consists of partial symmetric matrices 
with sparsity pattern $E$ that admit a positive semidefinite completion,
while the dual cone $K^*$ consists of positive semidefinite matrices
with the same sparsity pattern.
In the implementation, the pattern $E$ is obtained from a chordal extension
of the aggregate sparsity pattern of the problem, so that computations
with $K$ and~$K^*$ can be carried out using sparse Cholesky factorizations.

We use the logarithmic barrier kernel associated with $K^*$:
$\phi_*(S) = -\log\det(S)$.  The kernel on $K$ is the (negative) conjugate
\[
\phi(x) = \sup_{s\in\intr(K^*)} \{-s^Tx-\phi_*(s)\}.
\]
This choice is well adapted to sparse semidefinite cones because
its gradient and directional Hessian, as well as the corresponding Bregman
projection (\ref{e-num-proj}), can be evaluated efficiently
via sparse Cholesky factorization.

The kernel $\phi$ satisfies~(\ref{e-log-hom}) with $\nu=n$
and $\dom\phi=\intr(K)$, and is infinite on the boundary of~$K$.  
Theorem~\ref{thm-barrier} applies, but not Theorem~\ref{thm-finite},
and the exact-penalty objective converges at the rate 
$O((\log k)/k)$ for both methods.

\paragraph{Numerical results}
The sparse SDP instances are generated with the function \texttt{rand\_SDP}
from SMCP \cite{ADV:10}.  As in \cite{ADV:10}, the aggregate sparsity
patterns are taken from the SuiteSparse Matrix Collection \cite{KoA+:19}.
Chordal extensions are computed with CHOMPACK \cite{AnV:15}
using the AMD ordering \cite{ADD:96}.
CHOMPACK is also used for the sparse Cholesky factorizations
and related operations on the chordal embedding.
The remaining sparse matrix computations are implemented using CVXOPT
\cite{ADV:20}.

Table~\ref{t-sdp}%
\begin{table}
\centering
\begin{tabular}{rcccrrrrr} \toprule
  & & & & & & &
  \multicolumn{2}{c}{Projection statistics} \\
  \cmidrule(lr){8-9}
  Instance & $n$ & density & $m$ & Method & \# iter.\ &
  Time/iter.\ (s) & active & Newton \\ \midrule
  \multirow{2}{*}{rs35} & \multirow{2}{*}{$2003$}
    & \multirow{2}{*}{$20.9\%$} & \multirow{2}{*}{$2000$}
    & EG   & 21759 & 8.3520 & 16 & 5.36 \\
	& & & & PDHG & 22984 & 8.6134 & 27 & 5.29 \\
  \midrule
  \multirow{2}{*}{rs200} & \multirow{2}{*}{$3025$}
    & \multirow{2}{*}{$0.23\%$} & \multirow{2}{*}{$2000$}
    & EG   &  9578 & 0.7911 & 6 & 2.95 \\
	& & & & PDHG & 11038 & 0.9183 & 11 & 3.10 \\
  \midrule
  \multirow{2}{*}{rs228} & \multirow{2}{*}{$1919$}
    & \multirow{2}{*}{$0.88\%$} & \multirow{2}{*}{$2000$}
    & EG   &  5713 & 0.7622 & 8 & 3.04 \\
	& & & & PDHG & 8327 & 0.7513 & 7 & 2.99 \\
  \midrule
  \multirow{2}{*}{rs365} & \multirow{2}{*}{$4704$}
    & \multirow{2}{*}{$0.47\%$} & \multirow{2}{*}{$2000$}
    & EG   & 17474 & 3.7162 & 2 & 2.50 \\
	& & & & PDHG & 15037 & 3.8123 & 3 & 3.00 \\
  \midrule
  \multirow{2}{*}{rs1555} & \multirow{2}{*}{$7479$}
    & \multirow{2}{*}{$0.12\%$} & \multirow{2}{*}{$2000$}
    & EG   & 11759 & 1.4712 & 7 & 5.83 \\
	& & & & PDHG & 9187 & 1.5198 & 11 & 5.51 \\
  \bottomrule
\end{tabular}
\caption{Results for randomly generated sparse SDP instances.  The dimension
$n$ is the matrix order and $m$ is the number of equality constraints.  
The column \emph{active} gives the number of projections for which 
$\hat w^Tx=\mu$ at the solution (the total number of projections is twice 
\# iter.).  \emph{Newton} gives the average number of Newton steps
(averaged over the active projections). 
Inactive projections are available in closed form.}
\label{t-sdp}
\end{table}
reports the results for the sparse SDP instances.
As in Table~\ref{t-ntp}, the bound $\hat w^Tx\leq\mu$ is active at fewer
than 30 of the more than $10^4$ projections in each run, and those
projections take two to six Newton steps on average.
The main observation is that the cost per iteration is dominated
by the Cholesky factorizations and related sparse matrix operations  
in the Bregman projection oracle.  This keeps the
per-iteration cost inexpensive compared with methods based on dense
eigenvalue decomposition, and allows the first-order methods
to handle matrix variables with dimensions in the thousands.
Iteration counts of this size (about $10^4$ in Table~\ref{t-sdp}) are common
for first-order methods for SDPs,
even when only moderate accuracy is required.

\section{Conclusion} \label{s-conclusion}
We developed a bounded exact-penalty framework for applying Bregman
primal--dual proximal methods to conic linear programs.
The main idea is to replace the original conic LP by an equivalent
formulation over truncated cones.  Generalized projections on truncated
cones are often easier to define and compute than projections on the full
cones, and the truncation to bounded domains simplifies the convergence 
analysis.

For strictly feasible self-dual conic LPs, the required bounds are
available explicitly.  Combined with self-dual embeddings, this observation
gives a systematic procedure for obtaining bounded reformulations of general 
conic LPs, without
requiring problem-specific bounds on the primal or dual optimal sets.
The resulting problems can be solved by Bregman PDHG or Bregman extragradient
methods, with subproblems posed on truncated cones rather than on the full
cones.

We also analyzed the convergence of the Bregman PDHG and
extragradient methods with line search, applied to self-dual convex problems.
For self-dual conic LPs, the bounds are stated directly in terms of the 
exact-penalty objective.  They extend to distances generated by logarithmically homogeneous kernels, 
for which the Bregman radius on the truncated cone is infinite, 
and give an $O((\log k)/k)$ rate for the averaged iterates.
The rate improves to $O(1/k)$ if the kernel is finite at a solution
of the self-dual conic LP.

Numerical experiments on nonnegative trigonometric polynomial optimization
and sparse semidefinite optimization indicate that the approach is
promising for cone families with efficient Bregman projections.
Further algorithmic improvements are needed before it can serve
as a practical general-purpose first-order solver for structured conic LPs.
Two directions are restarting~\cite{AHLL:23,ADHL+:26,LXGY:25} and
preconditioning~\cite{PoC:11,ADHL+:21}.  Both have been highly effective
in first-order methods for LPs and conic LPs.
It remains to determine how these techniques should be adapted to
Bregman primal--dual proximal methods, and how they interact with the
exact-penalty formulations studied here.

\appendix

\section{Mixed-form conic LPs and exact penalties} \label{s-app-mixed}

This appendix records the mixed-form variants of the constructions used
in the main text.  The main text presents the exact-penalty reformulation
for conic LPs written without explicit free variables or equality blocks.
Here we give the corresponding formulation for mixed-form conic LPs, which
include free variables and explicit equality constraints.  This form is
also useful for self-dual LPs whose cones are not proper and for the
Ye--Todd--Mizuno embedding in Appendix~\ref{s-app-ytm}.

Consider the primal--dual pair
\begin{subequations}
\label{e-lp-mixed}
\BEA
\mbox{Primal:} &&
 \begin{array}[t]{ll}
 \mbox{minimize} & c_1^T x_1  + c_2^T x_2 \\
 \mbox{subject to} 
   & A_{11} x_1 + A_{12} x_2 - b_1 \in L_1 \\
   & A_{21} x_1 + A_{22} x_2 - b_2 = 0 \\
   & x_1 \in K_1,
\end{array}
\label{e-lp-mixed-primal}
\\
\mbox{Dual:} &&
 \begin{array}[t]{ll}
 \mbox{maximize} & b_1^T z_1 + b_2^T z_2 \\
 \mbox{subject to} 
   & -A_{11}^T z_1 - A_{21}^T z_2 + c_1 \in K_1^* \\
   & -A_{12}^T z_1 - A_{22}^T z_2 + c_2 = 0 \\
   & z_1 \in L_1^*.
\end{array}
\label{e-lp-mixed-dual}
\EEA
\end{subequations}
Here $K_1\subseteq\reals^{n_1}$ and $L_1\subseteq\reals^{m_1}$
are proper cones, $x_2\in\reals^{n_2}$ and $z_2\in\reals^{m_2}$ are
free variables, and the submatrices $A_{ij}$ have size $m_i\times n_j$.

\subsection{Strict consistency} \label{s-app-mixed-strict}

The primal problem~(\ref{e-lp-mixed-primal}) is a special case
of~(\ref{e-conic-lp-primal}) with
$K=K_1\times\reals^{n_2}$ and $L=L_1\times\{0\}$.  Hence
\[
K^* = (K_1 \times \reals^{n_2})^* = K_1^* \times \{0\}, 
\qquad 
L^*= (L_1 \times \{0\})^* = L_1^* \times \reals^{m_2},
\]
and 
\[
\lin(K^*) = \lin(K_1^*) \times \{0\} =  \{0\}\times \{0\}, \qquad
\lin(L^*) = \lin(L_1^*) \times \reals^{m_2}=  \{0\}\times \reals^{m_2}.
\]
The condition~(\ref{e-aff-hull-dim-m}) therefore means that the matrix
$\left[\begin{array}{cc} A_{21} & A_{22} \end{array}\right]$ has full row 
rank $m_2$.
The primal LP~(\ref{e-lp-mixed-primal}) is strictly consistent 
if and only if this rank condition holds and there exists
$\hat x=(\hat x_1,\hat x_2)\in\reals^{n_1}\times\reals^{n_2}$ such that
\[
 \hat x_1 \in\intr(K_1), \qquad
 A_{11}\hat x_1+A_{12}\hat x_2-b_1\in\intr(L_1), \qquad
 A_{21}\hat x_1+A_{22}\hat x_2-b_2=0.
\]
Similarly, the dual LP~(\ref{e-lp-mixed-dual}) is strictly consistent
if and only if
$\Rank\!\left(\begin{bmatrix}A_{12}^T&A_{22}^T\end{bmatrix}\right)=n_2$, and
there exists $\hat z=(\hat z_1,\hat z_2)\in\reals^{m_1}\times\reals^{m_2}$
such that
\[
 \hat z_1\in\intr(L_1^*), \qquad
 -A_{11}^T\hat z_1-A_{21}^T\hat z_2+c_1\in\intr(K_1^*), \qquad
 -A_{12}^T\hat z_1-A_{22}^T\hat z_2+c_2=0.
\]

\subsection{The self-dual case} \label{s-app-mixed-sd}

A useful special case of~(\ref{e-lp-mixed}) occurs when it is self-dual: 
$b=-c$, $L_1=K_1^*$, and $A = -A^T$.
That is, the problem (\ref{e-lp-mixed}) has the form
\[
\begin{array}{ll}
\mbox{minimize} &  -b_1^Tx_1 - b_2^Tx_2 \\
\mbox{subject to} 
& A_{11}x_1 + A_{12}x_2 - b_1 \in K_1^*, \\ 
& -A_{12}^Tx_1 + A_{22}x_2 - b_2 = 0, \\
& x_1 \in K_1,
\end{array}
\]
where $K_1$ is proper, $A_{11}=-A_{11}^T$, and
$A_{22}=-A_{22}^T$.

Splitting the free variable $x_2$ as $x_2=u_2-v_2$ gives the equivalent
self-dual LP over the proper cone
$K_1\times\reals^{n_2}_+\times\reals^{n_2}_+$:
\[
\begin{array}{ll}
\mbox{minimize} &  -b_1^Tx_1 - b_2^Tu_2 + b_2^T v_2 \\*[1ex]
\mbox{subject to} &
\left[\begin{array}{ccc}
 A_{11} & A_{12} & -A_{12} \\
-A_{12}^T & A_{22} & -A_{22} \\
 A_{12}^T & -A_{22} & A_{22}
\end{array}\right]
\left[\begin{array}{c} x_1 \\ u_2 \\ v_2 \end{array}\right]
-
\left[\begin{array}{c} b_1 \\ b_2 \\ -b_2 \end{array}\right]
\in K_1^* \times \reals^{n_2}_+ \times \reals^{n_2}_+,
\\*[3ex]
& (x_1,u_2,v_2)
\in K_1\times\reals^{n_2}_+\times\reals^{n_2}_+.
\end{array}
\]
The cone $K_1\times\reals^{n_2}_+\times\reals^{n_2}_+$ is proper.
Thus a self-dual LP over a nonproper cone can be represented as an
equivalent self-dual LP over a proper cone after making the lineality
space explicit.

\subsection{Exact penalties for mixed-form LPs}
\label{s-app-mixed-penalty}

The exact-penalty construction in Section~\ref{s-ex-pen} extends to the
mixed-form pair~(\ref{e-lp-mixed}).  The conic constraint blocks
are penalized by the corresponding Minkowski gauges, whereas the equality
blocks are penalized by dual norms.

Assume that the primal and dual optima are attained, and that there exist
optimal solutions
$x^\star=(x_1^\star,x_2^\star)$ and
$z^\star=(z_1^\star,z_2^\star)$ satisfying
\[
\hat w_1^T x_1^\star < 1, \qquad
\|x_2^\star\| < 1, \qquad
\hat s_1^T z_1^\star < 1, \qquad
\|z_2^\star\|' < 1,
\]
where $\|\cdot\|$ is a norm on $\reals^{n_2}$ and $\|\cdot\|'$ is a norm
on $\reals^{m_2}$.  Then the following penalized problems have the same
optimal values as the mixed-form primal and dual LPs.  Every optimal
solution of a penalized problem is optimal for the corresponding original
LP, and conversely every original optimal solution satisfying the imposed
bounds is optimal for the corresponding penalized problem:
\begin{subequations}
\label{e-mixed-lp-pen}
\BEA
\mbox{Primal:} && \begin{array}[t]{ll}
  \mbox{minimize} & c_1^T x_1 + c_2^Tx_2
	+ \sigma_{L_1}(b_1-A_{11}x_1-A_{12}x_2,\hat s_1)
	+ \|b_2-A_{21}x_1-A_{22}x_2\|_*' \\
  \mbox{subject to} & x_1\in K_1, \;\; \hat w_1^Tx_1\leq 1 \\
  & \|x_2\|\leq 1,
\end{array} \\
\mbox{Dual:} && \begin{array}[t]{ll}
  \mbox{maximize} & b_1^Tz_1+b_2^Tz_2
	-\sigma_{K_1^*}(A_{11}^Tz_1+A_{21}^Tz_2-c_1,\hat w_1)
	-\|{-A_{12}^Tz_1}-A_{22}^Tz_2+c_2\|_* \\
  \mbox{subject to} & z_1\in L_1^*, \;\; \hat s_1^Tz_1\leq 1 \\
  & \|z_2\|'\leq 1.
\end{array}
\EEA
\end{subequations}
In these problems, $\|\cdot\|_*$ and $\|\cdot\|_*'$ are the dual norms of
$\|\cdot\|$ and $\|\cdot\|'$.
The problems~(\ref{e-mixed-lp-pen}) are a special case of~(\ref{e-pd})
with $f=\delta_S$, the indicator function of the set
\[
S=\{(x_1,x_2)\in\reals^{n_1}\times\reals^{n_2}
\mid x_1\in K_1,\; \hat w_1^Tx_1\leq 1,\; \|x_2\|\leq 1\},
\]
and $g=\delta_T^*$, the support function of the set
\[
T=\{(z_1,z_2)\in\reals^{m_1}\times\reals^{m_2}
\mid z_1\in L_1^*,\; \hat s_1^Tz_1\leq 1,\; \|z_2\|'\leq 1\}.
\]

\section{Ye--Todd--Mizuno self-dual embedding} \label{s-app-ytm}
This appendix presents the Ye--Todd--Mizuno self-dual embedding for conic
LPs
\begin{subequations}
\label{e-primal-dual-gen}
\BEA 
\mbox{Primal:} & & 
\begin{array}[t]{ll}
\mbox{minimize} & c^T x \\
\mbox{subject to} & Ax - b = s \\
                  & x\in K, \;\; s \in L,
\end{array} \\
\mbox{Dual:} & & 
\begin{array}[t]{ll}
\mbox{maximize} & b^T z \\
\mbox{subject to} & -A^T z + c = w \\
  & z \in L^*, \;\; w \in K^*.
\end{array}
\EEA
\end{subequations}
The embedding converts the primal--dual pair into a single self-dual 
conic LP with known strictly feasible points and explicit bounds on the 
optimal set of the embedded problem. 

We assume the cones $K$ and $L$ are nonempty closed convex cones 
and $A\in\reals^{m\times n}$.
If the cones are not proper, one can make the linear hulls and 
lineality spaces of $K$ and $L$ explicit, which leads to a pair
of equality-constrained primal and dual problems of the form
\begin{subequations}
\label{e-primal-dual-gen-mixed}
\BEA 
\mbox{Primal:} & & 
\begin{array}[t]{ll}
\mbox{minimize} & c_1^T x_1 + c_2^T x_2 \\*[1ex]
\mbox{subject to} & 
\left[\begin{array}{cc} A_{11} & A_{12} \\ A_{21} & A_{22}
 \end{array}\right] \left[\begin{array}{c} x_1 \\ x_2 
\end{array}\right] - \left[\begin{array}{c} b_1 \\ b_2
 \end{array}\right] = 
 \left[\begin{array}{c} s_1 \\ s_2 \end{array}\right] \\*[2ex]
  & (x_1,x_2) \in K_1 \times \reals^{n_2}, \;\; 
  (s_1,s_2) \in L_1 \times \{0\},
\end{array} \\
\mbox{Dual:} & & 
\begin{array}[t]{ll}
\mbox{maximize} & b_1^T z_1 + b_2^Tz_2 \\*[1ex]
\mbox{subject to} & -\left[\begin{array}{cc} A_{11}^T & A_{21}^T \\
  A_{12}^T & A_{22}^T \end{array}\right] \left[\begin{array}{c}
 z_1 \\ z_2 \end{array} \right] + \left[\begin{array}{c} c_1 \\ c_2
 \end{array}\right] = \left[\begin{array}{c} w_1 \\ w_2 \end{array}
 \right]\\*[2ex]
  & (z_1,z_2) \in L_1^* \times \reals^{m_2}, \;\; 
    (w_1,w_2) \in K_1^* \times \{0\},
\end{array}
\EEA
\end{subequations}
where $K_1,L_1$ are proper cones in $\reals^{n_1}$ and $\reals^{m_1}$, 
respectively, and $A_{ij} \in\reals^{m_i\times n_j}$ for $i,j=1,2$.
Important special cases are LPs in standard form ($n_2=m_1=0$)
and LPs in inequality form ($n_1=m_2=0$).
It is useful to note that a problem in the 
form~(\ref{e-primal-dual-gen-mixed}) can be converted to the 
simpler~(\ref{e-primal-dual-gen}) with proper cones $K,L$, by the
familiar procedure of splitting free variables into nonnegative parts
and writing each equality as two opposing inequalities, for example, as
\begin{subequations}
\label{e-primal-dual-gen-mixed-std}
\BEA 
\mbox{Primal:} & & 
\begin{array}[t]{ll}
\mbox{minimize} & c_1^T x_1 + c_2^T x_{2+} - c_2^T x_{2-} \\*[1ex]
\mbox{subject to} & 
\left[\begin{array}{ccc} 
A_{11} & A_{12} & -A_{12} \\ 
A_{21} & A_{22} & -A_{22} \\
-A_{21} & -A_{22} & A_{22} 
 \end{array}\right] \left[\begin{array}{c} x_1 \\ x_{2+} \\ x_{2-}  
\end{array}\right] - \left[\begin{array}{c} b_1 \\ b_2 \\ -b_2 
\end{array}\right] = 
  \left[\begin{array}{c} s_1 \\ s_{2+} \\ s_{2-} \end{array}\right] \\*[2ex]
  & (x_1,x_{2+},x_{2-}) \in K_1 \times \reals_+^{n_2} \times \reals^{n_2}_+ \\
  & (s_1,s_{2+}, s_{2-}) \in L_1 \times \reals^{m_2}_+ \times \reals^{m_2}_+,
\end{array} \\
\mbox{Dual:} & & 
\begin{array}[t]{ll}
\mbox{maximize} & b_1^T z_1 + b_2^Tz_{2+} - b_2^T z_{2-} \\*[1ex]
\mbox{subject to} & -\left[\begin{array}{ccc} 
 A_{11}^T & A_{21}^T & -A_{21}^T \\
 A_{12}^T & A_{22}^T & -A_{22}^T \\
 -A_{12}^T & -A_{22}^T & A_{22}^T 
\end{array}\right] \left[\begin{array}{c}
 z_1 \\ z_{2+} \\ z_{2-} \end{array} \right] + 
\left[\begin{array}{c} c_1 \\ c_2 \\ -c_2 \end{array}\right] 
= \left[\begin{array}{c} w_1 \\ w_{2+} \\ w_{2-} \end{array} \right]\\*[2ex]
& (z_1,z_{2+}, z_{2-}) \in L_1^* \times \reals_+^{m_2} \times \reals_+^{m_2} \\
& (w_1,w_{2+}, w_{2-}) \in K_1^* \times \reals_+^{n_2} \times \reals_+^{n_2}.
\end{array}
\EEA
\end{subequations}

\subsection{The embedding} \label{s-app-ytm-embed}
The Ye--Todd--Mizuno self-dual embedding of~(\ref{e-primal-dual-gen}) is
\BEQ \label{e-sd-ext}
\begin{array}{ll}
\mbox{minimize} & \theta \\
\mbox{subject to} & 
 \left[\begin{array}{cccc}
0 & -A^T &  c  & q_x \\
A & 0    & -b & q_z \\
-c^T & b^T & 0 & q_{\tau} \\ 
-q_x^T & -q_z^T & -q_\tau & 0 
\end{array}\right]
\left[\begin{array}{c} x \\ z \\ \tau \\ \theta \end{array}\right]
+\left[\begin{array}{c} 0 \\ 0 \\ 0 \\ 1 \end{array}\right] = 
 \left[\begin{array}{c} w \\ s \\ \kappa \\ \eta \end{array}\right] 
 \\*[5ex]
& (x, z, \tau, \theta) \in K \times L^* \times \reals_+ \times \reals \\
& (w, s, \kappa, \eta) \in K^* \times L \times \reals_+ \times \{0\}.
\end{array}
\EEQ
The parameters $q_x,q_z,q_\tau$ are determined by the condition that 
given points 
\[
\hat x \in \relint(K), \qquad \hat z \in \relint(L^*), \qquad 
\hat w\in \relint(K^*), \qquad \hat s \in \relint(L), 
\]
are feasible with $\hat \tau =\hat\kappa = 1$, $\hat\eta =0$,
and $\hat\theta$ chosen so that the embedding equations are satisfied:
\[
\hat\theta = \hat w^T \hat x + \hat z^T \hat s + 1, \qquad
\left[\begin{array}{c} q_x \\ q_z \\ q_\tau \end{array}\right]
= \frac{1}{\hat\theta}
 \left( \left[\begin{array}{c} \hat w \\ \hat s \\ 1 \end{array}\right]
 - \left[\begin{array}{ccc} 0 & -A^T & c \\ A & 0 & -b \\
  -c^T & b^T & 0 \end{array}\right] 
 \left[\begin{array}{c} \hat x \\ \hat z \\ 1 \end{array}\right] \right).
\]
The embedding~(\ref{e-sd-ext}) is of the general form~(\ref{e-sd-lp}).
It can be written as
\[
 \begin{array}[t]{ll}
  \mbox{minimize} & -\tilde b^T \tilde x  \\
  \mbox{subject to} & \tilde A \tilde x - \tilde b = \tilde s \\
   & \tilde x\in \tilde K, \;\; \tilde s \in \tilde K^*,
 \end{array} 
\]
with 
\BEQ \label{e-ytm-AbK}
\tilde x = \left[\begin{array}{c} x \\ z \\ \tau \\ \theta
\end{array}\right], \qquad
\tilde s = \left[\begin{array}{c} w \\ s  \\ \kappa \\ \eta
\end{array}\right], \qquad
\tilde b = \left[\begin{array}{c} 0 \\ 0 \\ 0 \\ -1 \end{array}\right],
\qquad
\tilde A = \left[\begin{array}{cccc}
 0 & -A^T &  c & q_x \\
 A & 0    & -b & q_z \\
-c^T & b^T & 0 & q_{\tau} \\ 
-q_x^T & -q_z^T & -q_\tau & 0 
\end{array}\right]
\EEQ
and 
\BEQ \label{e-ytm-AbK-2}
\tilde K = K \times L^* \times \reals_+ \times \reals, 
\qquad
\tilde K^* = K^* \times L \times \reals_+ \times \{0\}.
\EEQ

The self-dual embedding of~(\ref{e-primal-dual-gen-mixed}) is
\BEQ \label{e-pd-gen-mixed-embedding}
\begin{array}{ll}
\mbox{minimize} & \theta \\
\mbox{subject to} & 
 \left[\begin{array}{cccccc}
0 & 0 & -A_{11}^T & -A_{21}^T &  c_1 & q_{x1} \\
0 & 0 & -A_{12}^T & -A_{22}^T &  c_2 & q_{x2} \\
A_{11} & A_{12} & 0  & 0  & -b_1 & q_{z1} \\
A_{21} & A_{22} & 0  & 0  & -b_2 & q_{z2} \\
-c_1^T & -c_2^T & b_1^T & b_2^T & 0 & q_{\tau} \\ 
-q_{x1}^T & -q_{x2}^T & -q_{z1}^T & -q_{z2}^T & -q_\tau & 0 
\end{array}\right]
\left[\begin{array}{c} x_1 \\ x_2 \\ z_1 \\ z_2 \\ \tau \\ \theta 
\end{array}\right]
+\left[\begin{array}{c} 0 \\ 0 \\ 0 \\ 0 \\ 0 \\ 1 
 \end{array}\right] = 
 \left[\begin{array}{c} w_1 \\ w_2 \\ s_1 \\ s_2 \\ \kappa \\ \eta 
 \end{array}\right] 
 \\*[8ex]
& (x_1, x_2, z_1, z_2, \tau, \theta) 
  \in K_1 \times \reals^{n_2} \times L_1^* \times \reals^{m_2} \times
 \reals_+ \times \reals \\
& (w_1, w_2, s_1, s_2, \kappa, \eta) \in 
  K_1^* \times \{0\} \times L_1 \times \{0\} \times \reals_+ \times  \{0\}.
\end{array}
\EEQ
The parameters $q_{x1}, q_{x2}, q_{z1}, q_{z2}, q_\tau$ are determined
by the condition that given points 
\[
(\hat x_1, \hat x_2) \in \intr(K_1) \times \reals^{n_2}, \quad
(\hat z_1, \hat z_2) \in \intr(L_1^*) \times \reals^{m_2}, \quad
(\hat w_1, \hat w_2) \in \intr(K_1^*) \times \{0\}, \quad
(\hat s_1, \hat s_2) \in \intr(L_1) \times \{0\}
\]
are feasible with $\hat\tau=\hat\kappa = 1$, $\hat \eta =0$, and
$\hat\theta = \hat w_1^T \hat x_1 + \hat s_1^T \hat z_1 + 1$. 
Here necessarily $\hat w_2=0$ and $\hat s_2=0$.
Problem~(\ref{e-pd-gen-mixed-embedding}) is also a self-dual problem of 
the form~(\ref{e-sd-lp}).

\paragraph{Strict consistency}
The self-dual embedding~(\ref{e-sd-ext}) satisfies the strong consistency 
condition~(\ref{e-sd-lp-strict-feas}) by construction, since
\BEAS
 (\hat x, \hat z, 1, \hat \theta) 
& \in & \relint(\tilde K) \\
& = & \relint(K) \times \relint(L^*) \times 
  \reals_{++} \times \reals, \\
(\hat w, \hat s, 1, \hat\eta) 
& \in & \relint(\tilde K^*) \\
& = & \relint(K^*) \times \relint(L) \times  \reals_{++} \times \{0\}
\EEAS
where $\tilde K$ is defined in~(\ref{e-ytm-AbK-2}).
The additional condition~(\ref{e-lin-K-K*}) needed for strict consistency, 
\BEQ
\label{e-lin-K-K*-2} 
 \left. \begin{array}{l}
 y\in \lin(\tilde K)  \\ \tilde A y \in \lin(\tilde K^*) 
 \end{array}\right\} \qquad \Longleftrightarrow \qquad y =0,
\EEQ
generally requires additional assumptions on $A, K, L$ 
in~(\ref{e-primal-dual-gen}).
\BIT
\item The embedding is strictly consistent if 
$K$ is pointed and $L$ has nonempty interior.  Then 
\[
\lin(\tilde K) = \{0\} \times \{0\} \times \{0\} \times \reals,
\qquad
\lin(\tilde K^*) = \lin(K^*) \times \lin(L) \times \{0\} \times \{0\},
\]
so the condition~(\ref{e-lin-K-K*-2}) requires that 
\[
\left[\begin{array}{cccc}
 0 & -A^T & c & q_x \\
 A &  0   & -b & q_z \\
-c^T & b^T & 0 & q_\tau \\
 -q_x^T & -q_z^T & -q_\tau & 0 
\end{array}\right] 
\left[\begin{array}{c} 0 \\ 0 \\ 0 \\ \theta  \end{array}\right]
 = \left[\begin{array}{c}  w \\ s \\ 0 \\ 0 \end{array}\right],
\qquad w\in \lin(K^*), \quad s \in \lin(L)
\]
holds only if $\theta = 0$.
To show this, we take the inner product with $(\hat x, \hat z, \hat\tau,
\hat\theta)$ on both sides of the equation:
\[
\hat x^T w + \hat z^T s 
= \left[\begin{array}{c} \hat x \\ \hat z \\ \hat \tau \\ \hat\theta
\end{array}\right]^T
\left[\begin{array}{cccc}
 0 & -A^T & c & q_x \\
 A &  0   & -b & q_z \\
-c^T & b^T & 0 & q_\tau \\
 -q_x^T & -q_z^T & -q_\tau & 0 
\end{array}\right] 
\left[\begin{array}{c} 0 \\ 0 \\ 0 \\ \theta  \end{array}\right]
= 
-\left[\begin{array}{c} \hat w \\ \hat s \\ \hat \kappa \\ -1
\end{array}\right]^T
\left[\begin{array}{c} 0 \\ 0 \\ 0 \\ \theta \end{array}\right]
= \theta.
\]
The two inner products on the left-hand side are zero
because $w\in\lin(K^*) = \span(K)^\perp$, $\hat x\in K$  and
$s \in \lin(L) = \span(L^*)^\perp$, $\hat z \in L^*$.

\item If $K$ is not pointed or $L$ has empty interior, 
the problem can be reformulated as the LP~(\ref{e-primal-dual-gen-mixed}),
where $K_1$ and $L_1$ are proper cones and the dimension of $A_{ij}$ is 
$m_i \times n_j$.  
The embedding~(\ref{e-pd-gen-mixed-embedding}) is strictly consistent if
\BEQ \label{e-strict-cons-ranks}
\Rank(\left[\begin{array}{cc} A_{21} & A_{22} \end{array}\right]) = m_2, 
\qquad
\Rank(\left[\begin{array}{c} A_{12} \\ A_{22} \end{array}\right]) = n_2. 
\EEQ
To see this, first note that the embedding is strongly consistent by 
construction, with known feasible points
\[
(\hat x_1, \hat x_2, \hat z_1, \hat z_2, \hat \tau, \hat \theta) 
\in \relint(K_1) \times \reals^{n_2} \times \relint(L_1^*) 
 \times \reals^{m_2} \times \reals_{++} \times \reals,
\]
and
\[
(\hat w_1, \hat w_2, \hat s_1, \hat s_2, \hat \kappa, \hat \eta) \in 
\relint(K_1^*) \times \{0\} \times \relint(L_1) \times \{0\} \times 
\reals_{++} \times \{0\},
\]
where $\hat\theta = \hat w_1^T \hat x_1 + \hat s_1^T \hat z_1 +1$.
The condition~(\ref{e-lin-K-K*-2}) requires that
\BEQ \label{e-gen-lp-Ay}
\left[\begin{array}{cccccc}
0 & 0 & -A_{11}^T & -A_{21}^T & c_1 & q_{x1} \\
0 & 0 & -A_{12}^T & -A_{22}^T & c_2 & q_{x2} \\
A_{11} & A_{12} & 0 & 0 & -b_1 & q_{z1} \\
A_{21} & A_{22} & 0 & 0 & -b_2 & q_{z2} \\
-c_1^T & -c_2^T & b_1^T & b_2^T & 0 & q_\tau \\
-q_{x1}^T & -q_{x2}^T & -q_{z1}^T & -q_{z2}^T & - q_\tau & 0
\end{array}\right] \left[\begin{array}{c}
 0 \\ x_2 \\ 0 \\ z_2 \\ 0 \\ \theta \end{array}\right] = 0 
\EEQ
implies $x_2=0$, $z_2=0$, $\theta = 0$.  
We first show that it implies $\theta =0$.
Taking the inner product of~(\ref{e-gen-lp-Ay})
with $(\hat x_1, \hat x_2, \hat z_1, \hat z_2, 1, \hat \theta)$,
we see that
\[
0= \left[\begin{array}{c} \hat w_1 \\  \hat w_2 \\ \hat s_1 \\
 \hat s_2 \\ 1 \\ -1 \end{array}\right]^T
 \left[\begin{array}{c}  0 \\ x_2 \\ 0 \\ z_2 \\ 0 \\ \theta \end{array}
 \right] 
= \hat w_2^T x_2 + \hat s_2^T z_2 -\theta = -\theta.
\]
If $\theta =0$ in~(\ref{e-gen-lp-Ay}), then
\[
\left[\begin{array}{cc}
 0 & -A_{21}^T \\
 0 & -A_{22}^T \\
A_{12} & 0 \\
A_{22} & 0 \\
-c_2^T & b_2^T \\
-q_{x2}^T & -q_{z2}^T \end{array}\right] \left[\begin{array}{c}
 x_2 \\ z_2 \end{array}\right] = 0,
\]
and the rank conditions~(\ref{e-strict-cons-ranks}) 
imply that $x_2=0$, $z_2=0$.

\item
The primal problem in~(\ref{e-primal-dual-gen}) 
is in the standard LP form if $L = \{0\}$ and $K$ is a proper cone.
The self-dual embedding is strictly consistent if $A$ has full row rank 
$m$.   
This is a special case of the rank conditions~(\ref{e-strict-cons-ranks}) 
with  $m_1=0$, $m_2=m$, $n_1 =n$, $n_2=0$.

\item 
The primal problem in~(\ref{e-primal-dual-gen})  is in inequality form
if $K = \reals^n$ and $L$ is a proper cone.
The self-dual embedding is strictly consistent if $A$ has rank $n$. 
This is a special case of the rank conditions~(\ref{e-strict-cons-ranks}) 
with  $m_1=m$, $m_2=0$, $n_1 =0$, $n_2=n$.
\EIT

\paragraph{Optimal solution}
The Ye--Todd--Mizuno embedding is a standard technique for 
initialization and infeasibility detection in implementations of
interior-point methods.
In this paper, it is used as a tool for constructing an equivalent conic 
LP with simple, known bounds on the optimal set.
With this in mind, it is useful to note that if the  primal and dual
LPs~(\ref{e-primal-dual-gen}) are feasible, then $\kappa$ must be 
zero at the optimum of~(\ref{e-sd-ext}).  If at the optimum $\kappa> 0$, 
then $\tau =0$ and the optimal solution $x,z$ of the embedded problem
satisfies
\[
 -A^T z \in\relint(K^*), \qquad Ax \in \relint(L), \qquad
 -c^T x  + b^Tz > 0, \qquad x \in K, \qquad z \in L^*. 
\]
If $b^Tz>0$ this shows that the primal problem in~(\ref{e-primal-dual-gen})
is infeasible.  If $c^Tx < 0$ it shows the dual problem is infeasible.
If $\tau>0$ at the optimum of~(\ref{e-sd-ext}), it is easily verified
that $x/\tau$, $s/\tau$, $w/\tau$, $z/\tau$ are optimal 
for~(\ref{e-primal-dual-gen}).

\subsection{Exact-penalty formulation}
If the embedding is strictly consistent, then its optimal value 
is zero and the optimal set is nonempty and bounded, as explained in 
Section~\ref{s-sd-oc}.  Optimal solutions $(x,z,\tau,\theta)$ and
$(w,s,\kappa,\eta)$ of~(\ref{e-sd-ext}) satisfy $\theta = 0$ and
\BEA
 \hat w^Tx + \hat s^T z + \tau 
 &  = & \hat w^Tx + \hat s^T z + \hat \kappa\tau  + \hat\eta\theta 
 \nonumber \\
 & = & 
\hat w^T\hat x + \hat s^T \hat z + \hat\kappa \hat\tau + \hat\eta \hat\theta
-w^T\hat x - s^T \hat z - \kappa\hat \tau - \eta\hat \theta 
 \nonumber \\
 & \leq & 
\hat w^T\hat x + \hat s^T \hat z + \hat \kappa \hat \tau + 
  \hat \eta\hat \theta \nonumber \\
 & = & \hat w^T\hat x + \hat s^T \hat z + 1 \nonumber \\
 & = & \hat\theta. \label{e-tym-bound}
\EEA
Hence, optimal $x,z,\tau$ satisfy
\[
x\in K, \qquad z \in L^*, \qquad 
\hat w^T x \leq \hat\theta, \qquad
\hat s^T z \leq \hat\theta, \qquad
0 \leq \tau \leq \hat \theta.
\]
Assume now that $K$ and $L^*$ are proper cones.  Then
$\hat w\in\intr(K^*)$ and $\hat s\in\intr(L)$.  For any
$\mu>\hat\theta$, the following exact-penalty formulation has the same
optimal value as~(\ref{e-sd-ext}) and preserves all embedded optimal
solutions satisfying the displayed bounds:
\[
\begin{array}{ll} 
\mbox{minimize} & \theta + 
\mu\sigma_{K^*}(A^Tz - \tau c - \theta q_x, \hat w)
+\mu\sigma_{L}(-Ax + \tau b - \theta q_z, \hat s)
\\ & \mbox{} 
+ \mu \max\{0, c^Tx - b^Tz -\theta q_\tau\}
+ |q_x^Tx + q_z^T z + q_\tau\tau -1| \\
\mbox{subject to} & 
 x\in K, \; \hat w^T x \leq \mu \\
 & z\in L^*, \; \hat s^T z \leq \mu \\
 & 0 \leq \tau \leq \mu \\
 & |\theta| \leq 1,
\end{array}
\] 
where $\mu > \hat\theta = \hat w^T \hat x + \hat z^T\hat s + 1$.

\section{Proofs of the convergence results} \label{s-app-proofs}
This appendix contains the proofs of the results in Section~\ref{s-ls}.

The proofs use the following property of Bregman proximal operators;
see, \eg, \cite{ChT:93,JiV:23}, \cite[property~1]{Tse:08}.
Let $y\in\intr(\dom\phi)$ and $a\in\reals^n$, and suppose
\[
 u = \argmin_x{(f(x) + a^Tx + \dd(x,y))}
\]
exists, is unique, and lies in $\intr(\dom\phi)$.  Then,
for all $x\in\dom f\cap\dom\phi$,
\BEQ \label{e-ls-prox-ineq}
 f(u) - f(x) + a^T(u-x) \leq \dd(x,y) - \dd(x,u) - \dd(u,y).
\EEQ

\subsection{Proof of Theorem~\ref{thm-ls-ergodic}} \label{s-thm-3}
We first derive an inequality for one iteration of
Algorithm~\ref{alg-ls}.  Fix $x,z\in D$ and $k\geq 0$.
The iterates $x_{k+1}$ and $z_{k+1}$ are the points~(\ref{e-ls-trial})
computed with $t=t_k$.
Applying~(\ref{e-ls-prox-ineq}) to~(\ref{e-ls-trial-x})
gives
\BEQ \label{e-ls-x-prox}
 t_k\left(f(x_{k+1}) - f(x) + (Az_k-b)^T(x_{k+1}-x)\right)
 \leq \dd(x,x_k) - \dd(x,x_{k+1}) - \dd(x_{k+1},x_k).
\EEQ
Similarly, applying~(\ref{e-ls-prox-ineq}) to~(\ref{e-ls-trial-z})
gives
\BEQ \label{e-ls-z-prox}
 t_k\left(f(z_{k+1}) - f(z) + (A(2x_{k+1}-x_k)-b)^T(z_{k+1}-z)\right)
 \leq \dd(z,z_k) - \dd(z,z_{k+1}) - \dd(z_{k+1},z_k).
\EEQ
Using the definition of $\mathcal L$ in~\eqref{e-lag-sd}
and $A^T = -A$, we write
\BEAS
\lefteqn{t_k\left(\mathcal L(x_{k+1},z) - \mathcal L(x,z_{k+1})\right)}\\
& = & t_k \big( f(x_{k+1}) - f(x) + f(z_{k+1}) - f(z)
  - b^T(x_{k+1}-x) - b^T(z_{k+1}-z) - z^TAx_{k+1} + z_{k+1}^TAx \big)\\
& = & t_k \big( f(x_{k+1}) - f(x) + (Az_k-b)^T(x_{k+1}-x)
  + f(z_{k+1}) - f(z) + (A(2x_{k+1}-x_k)-b)^T(z_{k+1}-z) \big) \\
& & \mbox{} + t_k (z-z_k)^TA(x-x_k) - t_k (z-z_{k+1})^TA(x-x_{k+1})
   - t_k (z_{k+1}-z_k)^TA(x_{k+1}-x_k) \\
& \leq &
 \dd(x,x_k) + \dd(z,z_k) + t_k (z-z_k)^TA(x-x_k)
 - \dd(x,x_{k+1}) - \dd(z,z_{k+1}) - t_k (z-z_{k+1})^TA(x-x_{k+1}) \\
& & \mbox{} - \dd(x_{k+1},x_k) - \dd(z_{k+1},z_k)
 - t_k (z_{k+1}-z_k)^TA(x_{k+1}-x_k).
\EEAS
The second equality can be verified by expanding the inner products
and using the skew-symmetry of $A$.  The inequality follows
from~(\ref{e-ls-x-prox}) and~(\ref{e-ls-z-prox}).
Finally, the line-search condition~(\ref{e-ls-test}), which holds
for the accepted stepsize $t_k$ and points $x_{k+1}, z_{k+1}$,
gives
\BEA
\lefteqn{t_k \left(\mathcal L(x_{k+1}, z) - \mathcal L(x,z_{k+1})\right)}
 \nonumber \\
& \leq &
 \dd(x,x_k) + \dd(z,z_k) + t_k (z-z_k)^TA(x-x_k)
 - \dd(x,x_{k+1}) - \dd(z,z_{k+1})
 - t_k (z-z_{k+1})^TA(x-x_{k+1})
 \nonumber \\
& & \mbox{} - (1-\delta)
 \left(\dd(x_{k+1},x_k) + \dd(z_{k+1},z_k)\right).
 \label{e-ls-one-step}
\EEA
Replacing $k$ by $i$ in~(\ref{e-ls-one-step}) and summing from $i=0$ to $i=k$,
and omitting the nonpositive terms with coefficient $-(1-\delta)$, gives
\BEA
\lefteqn{\sum_{i=0}^k t_i
 \left(\mathcal L(x_{i+1},z) - \mathcal L(x,z_{i+1})\right)} \nonumber \\
& \leq & \sum_{i=0}^k \big( \dd(x,x_i) + \dd(z,z_i)
 + t_i (z-z_i)^TA(x-x_i)
 - \dd(x,x_{i+1}) - \dd(z,z_{i+1}) - t_i (z-z_{i+1})^TA(x-x_{i+1}) \big)
 \nonumber \\
& = & \dd(x,x_0) + \dd(z,z_0) + t_0 (z-z_0)^TA(x-x_0)
 - \dd(x,x_{k+1}) - \dd(z,z_{k+1}) - t_k (z-z_{k+1})^TA(x-x_{k+1})
 \nonumber \\
& & \mbox{} + \sum_{i=1}^k (t_i - t_{i-1}) (z-z_i)^TA(x-x_i).
 \label{e-ls-conv-1}
\EEA
Since $x,z,x_i,z_i\in D$, each of the inner products
in~(\ref{e-ls-conv-1}) is bounded in absolute value by $\rho_2$.
Moreover,
\[
 t_0 + t_k + \sum_{i=1}^k |t_i - t_{i-1}|
 = 2 t_0 + 2\sum_{i=1}^k (t_i-t_{i-1})_+ = 2(t_0 + \Gamma_k),
\]
because $t_k - t_0 = \sum_{i=1}^k (t_i - t_{i-1})$.
Hence the inner product terms in~(\ref{e-ls-conv-1}) add up to at most
$2\rho_2(t_0+\Gamma_k)$.
Omitting the terms $-\dd(x,x_{k+1}) - \dd(z,z_{k+1}) \leq 0$, we obtain
\BEQ \label{e-ls-conv-2}
 \sum_{i=0}^k t_i
 \left(\mathcal L(x_{i+1},z) - \mathcal L(x,z_{i+1})\right)
 \leq \dd(x,x_0) + \dd(z,z_0) + 2\rho_2(t_0+\Gamma_k).
\EEQ
Since $\mathcal L(\cdot,z)$ is convex and $\mathcal L(x,\cdot)$
is concave, and $\xavg_k$, $\zavg_k$ are convex combinations of
$x_1,\ldots,x_{k+1}$ and $z_1,\ldots,z_{k+1}$ with weights
$t_i/T_k$,
\[
 \mathcal L(\xavg_k, z) - \mathcal L(x,\zavg_k) \leq
 \frac{1}{T_k} \sum_{i=0}^k t_i
 \left(\mathcal L(x_{i+1},z) - \mathcal L(x,z_{i+1})\right).
\]
Combined with~(\ref{e-ls-conv-2}) this gives~(\ref{e-ls-pointwise}).

The bound~(\ref{e-ls-gap-bound}) follows by taking
the supremum over $z\in D$ and the infimum over $x\in D$ 
on both sides of the inequality~(\ref{e-ls-pointwise}).

\subsection{Proof of Theorem~\ref{thm-eg-ergodic}} \label{s-thm-4}
Fix $x\in D$ and $i\geq 0$.  Applying (\ref{e-ls-prox-ineq}) to
the second update (\ref{e-extragradient-sd-alg2}) in
Algorithm~\ref{alg-eg-ls} at iteration $i$ (so that $t=t_i$,
$\bar x = \bar x_i$, and $x^+=x_{i+1}$), 
we obtain
\BEAS
 f(x) + (A\bar x_i - b)^Tx + \frac{1}{t_i} \dd(x,x_i)
 & \geq & f(x_{i+1}) + (A\bar x_i - b)^Tx_{i+1}
 + \frac{1}{t_i} \dd(x_{i+1},x_i) + \frac{1}{t_i}\dd(x,x_{i+1}) \\
 & \geq & f(\bar x_i) - b^T\bar x_i + \frac{1}{t_i}\dd(x,x_{i+1}).
\EEAS
The second inequality is the stepsize condition~(\ref{e-tseng-sd}).
Rearranging and using the definition of $\mathcal L$ in~\eqref{e-lag-sd}
gives
\[
 \mathcal L(\bar x_i, x)
 = f(\bar x_i) - f(x) + b^T(x-\bar x_i) - x^TA\bar x_i
 \leq \frac{1}{t_i} \left(\dd(x,x_i) - \dd(x,x_{i+1})\right).
\]
Multiplying by $t_i$ and summing from $i=0$ to $i=k$ gives
\[
 \sum_{i=0}^k t_i \mathcal L(\bar x_i,x) \leq
 \dd(x,x_0) - \dd(x,x_{k+1}) \leq \dd(x,x_0).
\]
Since $\mathcal L(\cdot,x)$ is convex and $\xavg_k$ is a convex
combination of $\bar x_0, \ldots, \bar x_k$ with weights $t_i/T_k$,
this implies
\[
 \mathcal L(\xavg_k, x) \leq \frac{\dd(x,x_0)}{T_k}
 \quad \mbox{for all $x\in D$}.
\]
Adding this inequality and the same inequality at $z\in D$, and using
$\mathcal L(x,\xavg_k) = -\mathcal L(\xavg_k,x)$,
gives~(\ref{e-eg-sd-ergodic}).

The bound~(\ref{e-eg-sd-gap}) follows by taking
the supremum over $z\in D$ and the infimum over $x\in D$
on both sides of the inequality~(\ref{e-eg-sd-ergodic}).

\subsection{Proof of Theorem~\ref{thm-barrier}} \label{s-thm-5}
We first note some consequences of~(\ref{e-log-hom}).   These properties
are well known for logarithmic barrier functions ($\dom \phi = \intr(K)$).
\BIT
\item $\nabla\phi(x)^T x = -\nu$ for all $x\in\intr(K)$.
This can be seen by differentiating~(\ref{e-log-hom}) with respect to 
$\lambda$ at $\lambda = 1$.
\item  $-\nabla\phi(x) \in K^*$ for all $x\in\intr(K)$.
This can be shown  as follows. If  $-\nabla \phi(x) \not\in K^*$,
there exists a $y\in \intr(K)$ with $\nabla \phi(x)^Ty > 0$.  
By convexity and logarithmic homogeneity,
\BEAS
\phi(ty) 
& \geq & \phi(x) + \nabla \phi(x)^T (ty-x) \\
& = & \phi(x) + t\nabla \phi(x)^T y + \nu \\
& > & \phi(x) + \nu,
\EEAS
which contradicts $\phi(ty) = \phi(y) - \nu\log t$ for all $t> 0$.

\item The identity in~(\ref{e-log-hom}) holds for all $x\in\dom\phi$:
if $x\in\dom\phi$ and $\lambda>0$, then $\lambda x\in\dom\phi$ and
$\phi(\lambda x) = \phi(x) - \nu\log\lambda$.
To see this, define $x_\epsilon = (1-\epsilon)x + \epsilon\hat x$ for
$\epsilon \in (0,1]$.
Since $x\in\dom\phi\subseteq K$ and $\hat x\in\intr(K)$,
the point $x_\epsilon$ is in $\intr(K)$, so~(\ref{e-log-hom}) gives
$\phi(\lambda x_\epsilon) = \phi(x_\epsilon) - \nu\log\lambda$.
The points $\hat x$ and $\lambda\hat x$ are in $\intr(K)$,
which is the interior of $\dom\phi$.
Since $\phi$ is closed, $\phi(x_\epsilon) \to \phi(x)$ and
$\phi(\lambda x_\epsilon) \to \phi(\lambda x)$ as
$\epsilon \to 0$ \cite[theorem~7.5]{Roc:70}.

\item $\phi$ is nonincreasing with respect to the cone order:
\BEQ \label{e-phi-monotone}
 \phi(v) \geq \phi(u)
\EEQ
if $u,v\in\intr(K)$ and $u-v\in K$, and also if $v\in\dom\phi$ and
$u-v\in\dom\phi$.
In the first case, the inequality follows from convexity and
$-\nabla\phi(u)\in K^*$:
$\phi(v) \geq \phi(u) + \nabla\phi(u)^T(v-u) \geq \phi(u)$.
In the second case, let $w = u-v$, and define
$\lambda_1 = \exp(-\phi(v)/\nu)$ and $\lambda_2 = \exp(-\phi(w)/\nu)$.
By the previous property, $\phi(v/\lambda_1) = \phi(w/\lambda_2) = 0$.
The point $u/(\lambda_1+\lambda_2)$ is a convex combination of
$v/\lambda_1$ and $w/\lambda_2$, so $\phi(u/(\lambda_1+\lambda_2))\leq 0$
by convexity.  Therefore,
\[
 \phi(u) = \phi\big(u/(\lambda_1+\lambda_2)\big)
 - \nu\log(\lambda_1+\lambda_2)
 \leq -\nu\log(\lambda_1+\lambda_2) < -\nu\log\lambda_1 = \phi(v).
\]
\EIT

\paragraph{Proof}
Fix $x,z\in C$ and $\xi\in(0,1]$, and define
$x_\xi = (1-\xi)x + \xi\hat x$ and $z_\xi = (1-\xi)z + \xi\hat x$.
Since $\hat x\in\intr(K)$ and $\xi>0$, these points are in
$C\cap\intr(K)\subseteq D$.
Since $x_\xi-\xi\hat x = (1-\xi)x\in K$, combining (\ref{e-phi-monotone})
and~(\ref{e-log-hom}) gives
\[
 \phi(x_\xi) \leq \phi(\xi\hat x) = \phi(\hat x) + \nu\log(1/\xi),
\]
and similarly $\phi(z_\xi)\leq\phi(\hat x)+\nu\log(1/\xi)$.
Since $x_\xi,z_\xi\in D$ and $\cl D = C$,
the definition~(\ref{e-rho3}) of $\rho_3$ gives
\BEAS
 \dd(x_\xi,x_0) + \dd(z_\xi,z_0)
 & = & \phi(x_\xi) + \phi(z_\xi) - \phi(x_0) - \phi(z_0)
 + \nabla\phi(x_0)^T(x_0-x_\xi) + \nabla\phi(z_0)^T(z_0-z_\xi) \\
 & \leq & 2\nu\log(1/\xi) + \rho_3.
\EEAS
Applying~(\ref{e-ls-pointwise}) at $(x_\xi,z_\xi)$ therefore gives
\BEQ \label{e-barrier-proof-1}
 \mathcal L(\xavg_k,z_\xi) - \mathcal L(x_\xi,\zavg_k)
 \leq \frac{2\nu\log(1/\xi) + \rho_3 + 2\rho_2(t_0+\Gamma_k)}{T_k},
\EEQ
where $\mathcal L$ is defined in~\eqref{e-L-sdlp}.
Since $\mathcal L$ is affine in each of its arguments on $C$,
\BEAS
 \mathcal L(\xavg_k,z) & = & 
\mathcal L(\xavg_k,z_\xi)
 + \xi\big(\mathcal L(\xavg_k, z) - \mathcal L(\xavg_k,\hat x)\big) \\
 & = & 
\mathcal L(\xavg_k,z_\xi) + \xi (z-\hat x)^T (b-A\xavg_k), \\
 -\mathcal L(x,\zavg_k) & = & -\mathcal L(x_\xi,\zavg_k)
 + \xi\big(\mathcal L(\hat x,\zavg_k) - \mathcal L(x,\zavg_k)\big)  \\
& = & -\mathcal L(x_\xi,\zavg_k) + \xi(x-\hat x)^T (b- A\zavg_k).
\EEAS
Using $A\hat x - b = \hat s$ and $A^T=-A$ (Assumption~\ref{as-B1}), 
we can write the last terms as
\BEAS
(z-\hat x)^T(b-A\xavg_k)
 & = & \hat s^T(\hat x-z) + (z-\hat x)^TA(\hat x-\xavg_k), \\
(x-\hat x)^T(b-A\zavg_k)
 & = &  \hat s^T(\hat x-x) + (x-\hat x)^TA(\hat x-\zavg_k).
\EEAS
Since $x,z\in K$ and $\hat s\in K^*$, the first terms on the
right-hand sides are at most $\hat s^T\hat x$.  The second terms are
at most $\rho_2$, since $\cl D = C$.  Hence
\BEAS
 \mathcal L(\xavg_k,z) & \leq & \mathcal L(\xavg_k,z_\xi)
 + \xi(\hat s^T\hat x + \rho_2), \\
 -\mathcal L(x,\zavg_k) & \leq & -\mathcal L(x_\xi,\zavg_k)
 + \xi(\hat s^T\hat x + \rho_2).
\EEAS
Adding the two inequalities and combining with~(\ref{e-barrier-proof-1}) 
gives
\[
 \mathcal L(\xavg_k,z) - \mathcal L(x,\zavg_k) \leq
 \frac{2\nu\log(1/\xi) + \rho_3 + 2\rho_2(t_0+\Gamma_k)}{T_k}
 + 2\xi(\hat s^T\hat x + \rho_2).
\]
The right-hand side does not depend on $x$ and $z$.  Taking the
supremum over $z\in C$ and the infimum over $x\in C$ on the
left-hand side, and using~(\ref{e-ls-gap-exact-pen}),
we obtain
\[
 q(\xavg_k) + q(\zavg_k) = \eta(\xavg_k,\zavg_k)
 \leq \frac{2\nu\log(1/\xi) + \rho_3 + 2\rho_2(t_0+\Gamma_k)}{T_k}
 + 2\xi(\hat s^T\hat x + \rho_2).
\]
The bound~(\ref{e-barrier-rate}) follows by taking
$\xi = 1/(k+2)$ and using $T_k \geq (k+1)t_{\min}$.

For Algorithm~\ref{alg-eg-ls}, the same argument applies with
$\zavg_k = \xavg_k$, $z_0 = x_0$, and~(\ref{e-eg-sd-ergodic}) in place
of~(\ref{e-ls-pointwise}).  Using~(\ref{e-ls-gap-exact-pen}),
we obtain
\[
 2 q(\xavg_k) = \eta(\xavg_k,\xavg_k) \leq 
\frac{2\nu\log(1/\xi) + \rho_3}{T_k} + 2\xi(\hat s^T\hat x + \rho_2).
\]
The bound~(\ref{e-barrier-rate-eg}) follows by taking
$\xi = 1/(k+2)$ and using $T_k\geq (k+1)t_{\min}$.

\subsection{Proof of Theorem~\ref{thm-finite}} \label{s-thm-6}
Fix $x,z\in D$, and define $\tilde x = (x+x^\star)/2$ and
$\tilde z = (z+x^\star)/2$.  These points are in $C$.
Since $x, x^\star\in\dom\phi$, we have $x/2, x^\star/2 \in \dom \phi$ and
$\phi(x^\star/2) = \phi(x^\star) + \nu\log 2$.
Inequality~(\ref{e-phi-monotone}), together with $\tilde x-x^\star/2=x/2$,
therefore gives
\[
 \phi(\tilde x) \leq \phi(x^\star/2) = \phi(x^\star) + \nu\log 2,
\]
and similarly $\phi(\tilde z)\leq\phi(x^\star)+\nu\log 2$.
Hence $\tilde x,\tilde z\in D$, and the definition~(\ref{e-rho4})
of $\rho_4$ gives
\[
 \dd(\tilde x,x_0) + \dd(\tilde z,z_0) \leq 2\nu\log 2 + \rho_4.
\]
Applying~(\ref{e-ls-pointwise}) at $(\tilde x,\tilde z)$ gives
\BEQ \label{e-finite-proof-1}
 \mathcal L(\xavg_k,\tilde z) - \mathcal L(\tilde x,\zavg_k)
 \leq \frac{2\nu\log 2 + \rho_4 + 2\rho_2(t_0+\Gamma_k)}{T_k}.
\EEQ
Since $\mathcal L$ is affine in each of its arguments on $C$,
\[
 \mathcal L(\xavg_k,\tilde z) = \frac{1}{2}\big(\mathcal L(\xavg_k,z)
 + \mathcal L(\xavg_k,x^\star)\big), \qquad
 \mathcal L(\tilde x,\zavg_k) = \frac{1}{2}\big(\mathcal L(x,\zavg_k)
 + \mathcal L(x^\star,\zavg_k)\big).
\]
The point $x^\star$ is optimal, so 
\[
 \mathcal L(x^\star,u)\leq
\sup_{v\in C}{\mathcal L(x^\star,v)} = q(x^\star) = 0
\]
for all $u\in C$.
In particular, $\mathcal L(x^\star,\zavg_k)\leq 0$ and
$\mathcal L(\xavg_k,x^\star) = -\mathcal L(x^\star,\xavg_k)\geq 0$.
Therefore,
\BEAS
 \mathcal L(\xavg_k,z) - \mathcal L(x,\zavg_k) 
 & \leq &
 2\big(\mathcal L(\xavg_k,\tilde z) - \mathcal L(\tilde x,\zavg_k)\big) \\
 & \leq & \frac{2\left(2\nu\log 2 + \rho_4 + 2\rho_2(t_0+\Gamma_k)\right)}{T_k}.
\EEAS
Taking the supremum over $z\in D$ and the infimum over $x\in D$, and
using~(\ref{e-ls-gap-exact-pen}) and $T_k\geq(k+1)t_{\min}$,
gives~(\ref{e-finite-rate}).

For Algorithm~\ref{alg-eg-ls}, the same argument applies with
$\zavg_k = \xavg_k$, $z_0 = x_0$, and~(\ref{e-eg-sd-ergodic}) in place
of~(\ref{e-ls-pointwise}), and gives~(\ref{e-finite-rate-eg}).

\end{document}